%% file: hmsvII0714.tex
\documentclass[12pt,letterpaper]{amsart}
\usepackage{epic,eepic,latexsym, amssymb, amscd, amsfonts, xypic, floatflt}

 \newlength{\baseunit}               
 \newcount{\numlines}                
\newcommand{\point}{\vspace{3mm}\par \noindent \refstepcounter{subsection}{\bf \thesubsection.} }
\newcommand{\tpoint}[1]{\vspace{3mm}\par \noindent \refstepcounter{subsection}{\bf \thesubsection.} 
  {\em #1. ---} }
\newcommand{\epoint}[1]{\vspace{3mm}\par \noindent \refstepcounter{subsection}{\bf \thesubsection.} 
  {\em #1.} }
\newcommand{\bpoint}[1]{\vspace{3mm}\par \noindent \refstepcounter{subsection}{\bf \thesubsection.} 
  {\bf #1.} }

\newcommand{\bpf}{\noindent {\em Proof.  }}
\newcommand{\epf}{\qed \vspace{+10pt}}

\newcommand{\SSS}{\mathfrak{S}}

\newcommand{\A}{\mathbb{A}}
\newcommand{\Z}{\mathbb{Z}}

\newcommand{\C}{\mathbb{C}}

\newcommand{\proj}{\mathbb P}

\newcommand{\Proj}{\operatorname{Proj}}
\newcommand{\oh}{{\mathcal{O}}}
\newcommand{\cm}{{\mathcal{M}}}
\newcommand{\cmbar}{\overline{\cm}}

\newcommand{\al}{\alpha}

\newcommand{\be}{\beta}
\newcommand{\ep}{\epsilon}
\newcommand{\ga}{\gamma}
\newcommand{\Ga}{\Gamma}

\newcommand{\De}{\Delta}
\newcommand{\si}{\sigma}
\newcommand{\la}{\lambda}
\newcommand{\La}{\Lambda}

\newcommand{\Pic}{\operatorname{Pic}}

\newcommand{\Sym}{\operatorname{Sym}}

\newcommand{\SP}{S}

\newcommand{\Segre}{\mathcal{S}_3}

\newcommand{\cited}{}

\newcommand{\vw}{\mathbf{w}}
\newcommand{\vdeg}{\operatorname{\mathbf{deg}}}

\newcommand{\lremind}[1]{{\bf[label:  #1]}}
\newcommand{\secretnote}[1]{}
\newcommand{\notation}[1]{}
\renewcommand{\lremind}[1]{{}}

\newcommand{\cut}[1]{}

\begin{document}
\pagestyle{plain}
\title{
The moduli space of $n$ points on the line is cut out by simple quadrics when
$n$ is not six
}
\author{Benjamin Howard, John Millson, Andrew Snowden and Ravi Vakil}
\thanks{The first two authors were partially supported by 
NSF grant DMS-0104006 and the last author
was partially supported by NSF Grant DMS-0228011.
}
\date{Friday, July 14, 2006.} 
\begin{abstract}
  A central question in invariant theory is that of
  determining the relations among invariants.  Geometric invariant
  theory quotients come with a natural ample line bundle, and hence
  often a natural projective embedding.  This question translates to
  determining the equations of the moduli space under this embedding.
  This note deals with one of the most classical quotients, the space
  of ordered points on the projective line.  We show that under
  any linearization, this quotient is cut out (scheme-theoretically)
  by a particularly simple set of quadric relations, with the single
  exception of the Segre cubic threefold (the space of six points with
  equal weight).  Unlike many facts in geometric invariant
  theory, these results (at least for the stable locus) are
  field-independent, and indeed work over the integers.
\end{abstract}
\maketitle
\tableofcontents

{\parskip=12pt 

\cut{
(1) On the third to last
 line I believe the equation should be
 so $X'_{\Gamma'}= X'_{\Gamma'}$
 The second $\Gamma$ should have a superscript ``prime'' ($\Gamma$ has been
 replaced by $\Gamma'$ and $\Gamma''$ and you are proving the
 required identity for each of these).

(2)  John July 13.
}

\section{Introduction}

We consider the space of $n$ (ordered) points on the projective line,
up to automorphisms of the line.  In characteristic $0$, the best
description  of this is as the Geometric Invariant Theory quotient
$(\proj^1)^n // SL(2)$.  This is one of the most classical
examples of a GIT quotient, and is one of the first examples given in
any course (see \cite[\S 2]{ms}, \cite[\S 3]{git}, \cite[\S 4.5]{n},
\cite[Ch.~11]{d}, \cite[Ch.~I]{do}, ...). 
We find it remarkable that well over a century after this fundamental
invariant theory problem first arose, we have so little understanding
of its defining equations, when they should be expected to be quite simple.  

The construction of the quotient depends on a linearization, in the
form of a choice of weights of the points $\vw = (w_1, \dots, w_n)$
(the {\em weight vector}).  We denote the resulting quotient
$M_{\vw}$.  A particularly interesting case is when all points are
treated equally, when $\vw = (1, \dots, 1) = 1^n$.  We call this the
{\em equilateral} case because of its interpretation as the moduli
space of equilateral space polygons in symplectic geometry
(e.g.\ \cite[\S 2.3]{hmsv}). At risk of confusion, we denote this important
case $M_{1^n}$ by $M_n$ for simplicity.  We say that $n$
points of $\proj^1$ are {\em $\vw$-stable} (resp.\ {\em
  $\vw$-semistable}) if the sum of the weights of points that coincide
is less than (resp.\ no more than) half the total weight.  The
dependence on $\vw$ will be clear from the context, so the prefix
$\vw$- will be omitted.  The $n$ points are {\em strictly semistable}
if they are semistable but not stable.  Then $M_{\vw}$ is a projective
variety, and GIT gives a natural projective embedding.  The stable
locus of $M_{\vw}$ is a fine moduli space for the stable points of
$(\proj^1)^n$.  The strictly semistable locus of $M_{\vw}$ is a finite
set of points, which are the only singular points of $M_{\vw}$.  The
question we wish to address is: what are the equations of $M_{\vw}$?

We prefer to work as generally as possible, over the integers, so we
now define the moduli problem of stable $n$-tuples of points in
$\proj^1$.  For any scheme $B$, {\em families of stable $n$-tuples of
  points in $\proj^1$ over $B$} are defined to be $n$ morphisms to
$\proj^1$ such that no more than half the weight is concentrated at a
single point of $\proj^1$.  More precisely, a family is a morphism
$(\phi_1, \dots, \phi_n): B \times \{ 1, \dots, n \} \rightarrow
\proj^1$ such that for any $I \subset \{ 1, \dots, n \}$ such that
$\sum_{i \in I} w_i \geq \sum_{i=1}^n w_i / 2$, we have $\cap_{i \in
  I} \phi^{-1}_i(p) = \emptyset$ for all $p$.  Then there is a fine
moduli space for this moduli problem, quasiprojective over $\Z$, which
indeed has a natural ample line bundle, the one suggested by GIT.
(This is well-known, but in any case will fall out of our analysis.)
In parallel with GIT, we define stable points.  Hence this space has
extrinsic projective geometry.  The question we will address in this
context is: what are its equations?

The main moral of this note is taken from Chevalley's construction of
Chevalley groups: when understanding a vector space defined
geometrically, choosing a basis may obscure its structure.  Instead,
it is better to work with an equivariant generating set, and
equivariant linear relations. A prototypical example is the standard
representation of $\SSS_n$, which is best understood as the
permutation representation on the vector space generated by $e_1$,
\dots, $e_n$ subject to the relation $e_1 + \cdots + e_n = 0$.  As an
example, we give yet another short proof of Kempe's theorem (Theorem
~\ref{kempe}), which has been called the ``deepest result'' of
classical invariant theory \cite[p.~156]{howe} (in the sense that it
is the only result that Howe could not prove from standard
constructions in representation theory).  Another application is the
computation in \S \ref{degree} of the degree of all $M_{\vw}$.

We now state our main theorem.  We will describe a natural
(equivariant) set of ``graphical'' generators of the algebra of
invariants (in Section~\ref{snooker}).  The algebraic structure of the
invariants is particularly transparent in this language, and as an
example we give a short proof of Kempe's Theorem~\ref{kempe}, and give
an easy basis of the group of invariants (by ``non-crossing variables'', 
Theorem~\ref{basis}).  We
then describe some geometrically (or combinatorially) obvious
relations, the (linear) {\em sign relations}, the (linear) {\em
  Pl\"ucker relations}, and the (quadratic) {\em simple binomial
  relations}. 

\tpoint{Main Theorem} {\em 
With the single exception of $\vw =
(1,1,1,1,1,1)$, the following holds.
\begin{enumerate}\item[(a)]
Over a field of characteristic $0$, the GIT
  quotient $(\proj^1)^n // SL_2$ (with its natural projective embedding)
  is cut out scheme-theoretically by the sign, Pl\"ucker, and simple binomial
  relations.  
\item [(b)] Over $\Z$ (or any base scheme), the fine moduli space of stable
  $n$-tuples of points on $\proj^1$ is quasiprojective over $\Z$.
  Under its natural embedding, its closure is cut out by the sign,
  Pl\"ucker, and simple binomial relations.
\end{enumerate}
}
\label{mainthm}\lremind{mainthm}

The exceptional case $\vw = (1,1,1,1,1,1)$ is the Segre
cubic threefold. 

Note that Geometric Invariant Theory  does not apply to $SL(2)$-quotients
in positive characteristic, as $SL(2)$ is not a reductive group
in that case.  Thus  for (b) we must construct the moduli space
by other means.  

The idea of the proof is as follows.  We first reduce the question to
the equilateral case, where $n$ is even.  We do this by showing a stronger
result, which reduces such questions about the {\em ideal of relations}
of invariants to the equilateral case.

\tpoint{Theorem (reduction to equilateral case, informal statement)}
{\em For any weight $\vw$, there is a natural map from 
the graded ring of projective invariants for $1^{ | \vw |}$ to those
for $\vw$.  Under this map, each of our generators for
$1^{|\vw|}$ is sent to either a generator for $\vw$, or to zero.
Moreover, the relations for $1^{| \vw|}$ generate the relations
for $\vw$.
\label{preliminary}\lremind{preliminary}
}

In Section \ref{undem}, we will state this result precisely, and prove
it, once we have introduced some terminology.  This is a stronger result than we will need, as it refers to the
full ring of projective invariants, not just to the quotient variety.
It is also a stronger statement than simply saying that $M_{\vw}$ is
naturally a linear section of $M_{| \vw|}$.

We then  verify
Theorem~\ref{mainthm} ``by hand'' in the cases $n=2m \leq 8$ (\S
\ref{smallcases}).  
The cases $n=6$ and $n=8$ are the base cases for
our later argument (which is ironic, as in the $n=6$ case the result
does not hold!).  In \S \ref{proof}, we show that the result holds
set-theoretically, and that the projective variety is a fine moduli
space away from the strictly semistable points.  The strictly
semistable points are more delicate, as the quotient is not naturally
a fine moduli space there; we instead give an explicit description of
a neighborhood of a strictly semistable point, as the affine variety
corresponding to rank one $(n/2-1) \times (n/2-1)$ matrices with
entries distinct from $1$, using the Gel'fand-MacPherson
correspondence.  We prove the result in this neighborhood in \S
\ref{sss}.

\bpoint{Related questions}
Our question is related to another central problem in invariant
theory: the invariants of binary forms, or equivalently $n$ unordered
points on $\proj^1$, or equivalently equations for $M_n / \SSS_n$.  These
generators and relations are more difficult, and the relations are certainly not just
quadratic.  Mumford describes Shioda's solution for $n=8$
\cite{shioda} as ``an extraordinary tour de force'' \cite[p.~77]{git}.
One might dream that the case $n=10$ might be tractable by computer,
given the explicit relations for $M_{10}$ described here.

We remark on the relation of this paper to our previous paper
\cite{hmsv}.  That paper dealt with the {\em ideal} of relations among
the invariants, and a central result was that this ideal was cut out
by equations of degree at most $4$.  The present article deals with
the moduli space as a projective variety; it is not at all clear that
our quadrics generate the ideal of invariants.

{\bf Question.}  Do the simple binomial relations generate the ideal of
relations among invariants (if $n \neq 6$)?

By Theorem~\ref{preliminary}, it suffices to consider the equilateral
case.
\cite{hmsv} suggests one approach: one could hope to show that the
explicit generators given there lie in the ideal given by these
simple binomial quadrics.

Even special cases are striking, and are simple to state but
computationally too complex to verify even by computer.
For example, we will describe particularly attractive
relations for $M_{10}$ of degree
degrees $3$, $5$, $7$ (\S \ref{sneaky2}); do these
lie in the ideal of our simple quadrics?  We
will describe a quadric relation for $M_{12}$ (\S \ref{sneaky1}).
Is this a linear combination of our simple binomial quadrics? 

These questions lead to more speculation.  The shape of the proof of
the main theorem suggests an explanation for the existence of an
exception: when the number of points gets ``large enough'', the
relations are all inherited from ``smaller'' moduli spaces, so at some
point (in our case, when $n=8$) the relations ``stabilize''.

{\bf Speculative question.}  We know that $M_n$ satisfies
Green's property $N_0$ (projectively normality --- Kempe's Theorem~\ref{kempe}),
and the questions above suggest that $M_n$ might satisfy property
$N_1$ (projectively normal and cut out by quadrics) for $n >6$.  Might it be
true that each $p$, $M_n$ satisfies $N_p$ for $p \gg 0$?
(Property $N_2$ means that the scheme satisfies $N_1$, and the
syzygies among the quadrics are linear.  These properties measure
the ``niceness'' of the ideal of relations.)

{\bf Speculative question.}  Might similar results hold true for other
moduli spaces of a similar flavor, e.g.\ $\cmbar_{0,n}$ or even
$\cmbar_{g,n}$?  Are equations for moduli spaces inherited from
equations of smaller moduli spaces for $n \gg 0$?  For further
motivation, see the striking work of \cite{kt} on equations
cutting out $\cmbar_{0,n}$.  
Also, an earlier example of quadratic equations inherited
from smaller moduli spaces appeared in \cite{bp}, in the
case of Cox rings of del Pezzo surfaces.  I thank B.~Hassett
for pointing out this reference.

\noindent {\bf Acknowledgments.}  The last author thanks Lucia
Caporaso, Igor Dolgachev, Diane Maclagan, Brendan Hassett, and Harm
Derksen for useful suggestions.

\section{The invariants of $n$ points on $\proj^1$ as a graphical
  algebra}
\label{snooker}\lremind{snooker}
We give a convenient alternate description of the
generators (as a group) of the ring of invariants of $n$ ordered
points on $\proj^1$.  By {\em graph} we
will mean a {\em directed} graph on $n$ vertices labeled $1$ through
$n$.  Graphs may  have multiple edges, but may have no
loops.  The {\em multidegree} of a graph $\Ga$ is the $n$-tuple of
valences of the graph, denoted $\vdeg \Ga$.  The bold font is a
reminder that this is a vector.  We consider each graph as a set of
edges. For each edge $e$ of $\Ga$, let $h(e)$ be the head vertex of
$e$ and $t(e)$ be the tail.  We use multiplicative notation for the
``union'' of two graphs: if $\Ga$ and $\De$ are two graphs on the same
set of vertices, the union is denoted by $\Ga \cdot \De$ 
(so for example $\vdeg \Ga + \vdeg \De =
\vdeg \Ga \cdot \De$), see Figure~\ref{snip}.  (We will occasionally
use additive and subtractive notation when we wish to ``subtract''
graphs.  We apologize for this awkwardness.)

\begin{figure}
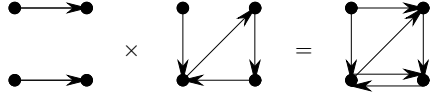

\begin{center}
\include{snip}
\end{center}
\caption{Multiplying (directed) graphs \lremind{snip}\label{snip}}
\end{figure}

For each graph $\Ga$, define
$X_\Ga \in H^0((\proj^1)^n, \oh_{(\proj^1)^n}(\vdeg \Ga))$
by \lremind{snack}
\begin{equation}X_\Ga = \prod_{\text{edge $e$ of $\Ga$}} (p_{h(e)} - p_{t(e)})
= \prod_{\text{edge $e$ of $\Ga$}} (u_{h(e)} v_{t(e)}- u_{t(e)} v_{h(e)}).
\label{snack}
\end{equation}
If $S$ is a non-empty set of graphs of
the same degree, $[X_{\Ga}]_{\Ga \in S}$ denotes a point in projective
space $\proj^{ \left| S \right| -1}$ (assuming some such $X_{\Ga}$ is
nonzero of course).  For any such $S$, the map $(\proj^1)^n
\dashrightarrow [X_{\Ga}]_{\Ga \in S}$ is easily seen to be invariant
under $SL(2)$: replacing $p_i$ by $p_i + a$ preserves
$X_{\Ga}$; replacing $p_i$ by $a p_i$ changes each $X_{\Ga}$ by the
same factor; and replacing $p_i$ by $1/p_i$ also changes each
$X_{\Ga}$ by the same factor.

The First Fundamental Theorem of Invariant Theory \cite[Thm.~2.1]{d}
states that, given a weight $\vw$, the ring of invariants of
$(\proj^1)^n // SL_2$ is generated (as a group) by the $X_{\Ga}$ where
$\vdeg \Ga$ is a multiple of $\vw$.  The translation to the tableaux
is as follows. Choose any ordering of the
edges $e_1$, \dots, $e_{\left| \Ga \right|}$ of $\Ga$.  Then $X_{\Ga}$ corresponds
to any $2 \times \left| \Ga \right|$ tableau where the top row of the $i$th
column is $h(e_i)$ and the bottom row is $t(e_i)$.  We will soon see
advantages of this graphical description as compared to the tableaux
description.

We now describe several types of relations among the $X_{\Ga}$, which
will all be straightforward: the sign relations, the Pl\"ucker (or
straightening) relations, the simple binomial relations, and
the Segre cubic relation.  

\epoint{The sign (linear) relations}
The {\em sign relation} $X_{\Ga \cdot \vec{xy}} = - X_{\Ga \cdot \vec{yx}}$
(Figure~\ref{snoopy})
is immediate, given the definition \eqref{snack}. 
Because of the sign relation, we may omit arrowheads
in identities where it is clear how to consistently insert them
(see for example Figures~\ref{s5} and \ref{s835}, where 
even the vertices are implicit).

\begin{figure}
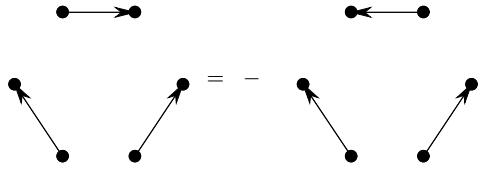

\begin{center}
\include{snoopy}
\end{center}
\caption{An example of the sign relation \lremind{snoopy}\label{snoopy}
}
\end{figure}

\epoint{The Pl\"ucker (linear) relations}
The identity of Figure~\ref{snide} may be verified by direct
calculation. 
If $\Ga$ is any graph on $n$ vertices, and
$\De_1$, $\De_2$, $\De_3$ are three graphs on the same vertices given by
identifying the four vertices of Figure~\ref{snide} with the some four
of the $n$ vertices of $\Ga$, then  \lremind{plucker}
\begin{equation}X_{\Ga \cdot \De_1} + X_{\Ga \cdot \De_2} + X_{\Ga \cdot \De_3} =
0.\label{plucker}
\end{equation}
These relations are called {\em Pl\"ucker relations} (or {\em
  straightening rules}).  See Figure~\ref{snoopdog} for an example.
We will sometimes refer to this relation as the Pl\"ucker relation for
$\Ga \cdot \De_1$ with respect to the vertices of $\De_1$.  

\begin{figure}
\begin{center}
\include{snide}
\end{center}
\caption{The Pl\"ucker relation for $n=4$ (and $\vw = (1,1,1,1)$)
\label{snide}\lremind{snide}}
\end{figure}

\begin{figure}
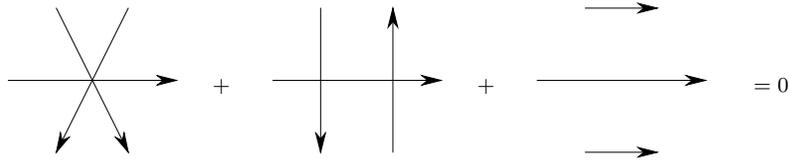

\begin{center}
\include{snoopdog}
\end{center}
\caption{An example of a Pl\"ucker relation 
\label{snoopdog}\lremind{snoopdog}}
\end{figure}

Using the Pl\"ucker relations, one can reduce the number of generators
to a smaller set, which we will do shortly (Proposition~\ref{snap}).
However, a central thesis of this article is that this is the wrong thing
to do too soon; not only does it obscure the $\SSS_n$ symmetry of this
generating set, it also makes certain facts opaque.  As an example, we
give a new proof of Kempe's theorem.  The proof will also serve as preparation for the
proof of the main theorem, Theorem~\ref{mainthm}.

\tpoint{Kempe's Theorem \cite[Thm.~4.6]{hmsv}} {\em 
The lowest degree invariants generate the ring of invariants.}
\lremind{kempe}\label{kempe}

Note that the lowest-degree invariants are of weight $\epsilon_{\vw}
\vw$, where $\ep_{\vw} = 1$ if $\left| \vw \right|$ is even, and $2$
if $\left| \vw \right|$ is odd.

\bpf We begin in the case when $\vw = (1, \dots, 1)$ where $n$ is even.  Recall Hall's Marriage
Theorem: given a finite set of men $M$ and women $W$, and some men and
women are compatible (a subset of $M \times W$), and it is desired to
pair the women and men compatibly, then it is necessary and sufficient
that for each subset $S$ of women, the number of men compatible with
at least one of them is at least $\left| S \right|$.  

Given a graph $\Ga$ of multidegree $(d, \dots, d)$, we show that we
can find an expression $\Ga = \sum \pm \De_i \cdot \Xi_i$ where $\vdeg
\De_i = (1, \dots, 1)$.  Divide the vertices into two equal-sized
sets, one called the ``positive'' vertices and one called the
``negative'' vertices.  This creates three types of edges: positive
edges (both vertices positive), negative edges (both vertices
negative), and neutral edges (one vertex of each sort).  When one applies the
Pl\"ucker relation to a positive edge and a negative edge, all
resulting edges are neutral (see Figure~\ref{snide}, and take two of
the vertices to be of each type).  Also, each regular graph must have
the same number of positive and negative edges.  
Working inductively on the number of positive edges, 
we can  use the
Pl\"ucker relations so that all resulting graphs have only neutral
edges. 
We thus have an expression $\Ga = \sum \pm \Ga_i$ where each 
$\Ga_i$ has only neutral edges and is hence a bipartite graph.
Each vertex of $\Ga_i$ has the same valence
$d$, so any set of $p$ positive vertices must connect to at least $p$ negative
edges. By
Hall's Marriage Theorem, we can find a matching $\De_i$ that is a subgraph
of $\Ga_i$, with ``residual
graph'' $\Xi_i$ (i.e.\, $\Ga_i = \De_i \cdot \Xi_i$).  Thus
the result holds in the equilateral case.



We next treat the general case. If $\left| \vw \right|$ is odd, it
suffices to consider the case $2 \vw$, so by replacing $\vw$ by $2
\vw$ if necessary, we may assume $\ep_{\vw}=1$.  The key idea is that
$M_{\vw}$ is a linear section of $M_{\left| \vw \right|}$.  Suppose
$\deg \Ga = d \vw$.  Construct an auxiliary graph $\Ga'$ on $\left|
\vw \right|$ vertices, and  a
map of graphs $\pi: \Ga' \rightarrow \Ga$ such that (i) the preimage of
vertex $i$ of $\Ga$ consists of $w_i$ vertices of $\Ga'$, (ii) $\pi$
gives a bijection of edges, and (iii) each vertex of $\Ga'$ has
valence $d$, i.e.\ $\Ga'$ is $d$-regular.  (See Figure~\ref{snifter}
for an illustrative example. There may be choice in defining $\Ga'$).
Then apply the algorithm of the previous paragraph to $\Ga'$.  By
taking the image under $\pi$, we have our desired result for $\Ga$.
\epf

\begin{figure}
\begin{center}
\include{snifter}
\end{center}
\caption{Constructing $\Ga'$ from $\Ga$ (example with
$\vw = (1,1,2,2)$, $d=2$) \label{snifter}\lremind{snifter}
}
\end{figure}

Choosing a planar representation makes termination of certain
algorithms straightforward as well, as illustrated by the following
argument.  Consider the vertices of the graph to be the vertices of a
regular $n$-gon, numbered (clockwise) $1$ through $n$.  A graph is
said to be {\em non-crossing} if no two edges cross.  Two edges
sharing one or two vertices are considered not to cross.
A variable $X_{\Ga}$ is said to be non-crossing if $\Ga$ is.

\tpoint{Proposition (graphical version of ``straightening algorithm'')}
{\em For each $\vw$, the non-crossing variables of degree $\vw$ 
generate $\langle X_{\Ga} \rangle_{\vdeg \Ga = \vw}$ (as an abelian group).}
\lremind{snap}
\label{snap}

This is essentially the straightening algorithm (e.g.\ \cite[\S
2.4]{d}) in this situation.

\bpf We explain how to express $X_{\Ga}$ in terms of non-crossing
variables.  If $\Ga$ has a crossing, choose one crossing $wx \cdot yz$
(say $\Ga = wx \cdot yz \cdot \Ga'$), and use the Pl\"ucker relation
\eqref{plucker} involving $wxyz$ to express $\Ga$ in terms of two other graphs $wy
\cdot xz \cdot \Ga'$ and $wz \cdot xy \cdot \Ga'$.  Repeat this if
possible.  We now show that this process terminates, i.e.\ that this
algorithm will express $X_{\Ga}$ in terms of non-crossing variables.
Both of these graphs have lower sum of edge-lengths than $\Ga$ (see
Figure~\ref{sniffle}, using the triangle inequality on the two
triangles with side lengths $a$, $d$, $f$ and $b$, $c$, $e$).  As
there are finite number of graphs of weight $\vw$, and hence a finite
number of possible sums of edge-lengths, the process must terminate.
\epf

\begin{figure}
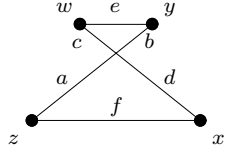

\begin{center}
\include{sniffle}
\end{center}
\caption{The triangle inequality implies termination of straightening:
$b+c>e$, $a+d>f$
\label{sniffle}\lremind{sniffle}}
\end{figure}

\tpoint{Theorem (non-crossing basis of invariants)}
{\em For each $\vw$, the non-crossing variables of degree $\vw$ 
form a basis for  $\langle X_{\Ga} \rangle_{\vdeg \Ga = \vw}$.}
\lremind{basis}
\label{basis}

\bpf Proposition~\ref{snap} shows that the non-crossing variables
span, so it remains to show that they are
linearly independent.  Assume otherwise that they are not always
linearly independent, and thus there s a simplest nontrivial relation
$R$ for some smallest $\vw$ (where the $\vw$ are partially ordered by
$\left| \vw \right|$ and $\# \vw$).  The relation $R$ states that some
linear combination of these graphs is $0$.  Suppose this relation involves
$\# \vw = n$ vertices.  Then not all of the graphs in $R$ contain edge
$(n-1)n$ (or else we could remove one copy of the edge from each term in $R$ 
and get a smaller
relation).

Then identify vertices $(n-1)$ and $n$, throwing out the graphs containing
edge $(n-1)n$.  This ``contraction'' gives a natural
bijection between non-crossing graphs on $n$ vertices with multidegree
$\vw$, not containing edge $(n-1)n$, and non-crossing graphs
on $n-1$ vertices with multidegree $\vw'$ (where $w'_{n-1} = w_{n-1}+w_n$).
The resulting relation still holds true (in terms of invariants,
we insert the relation $p_{n-1}=p_n$ into the older relation).  
We have contradicted the minimality of $R$, so no such 
counterexample $R$ exists.
\epf


\epoint{Binomial (quadratic) relations} We next describe some obvious
binomial relations.  If $\vdeg \Ga_1 = \vdeg \Ga_2$ and $\vdeg \De_1 =
\vdeg \De_2$, then clearly $X_{\Ga_1 \cdot \De_1 } X_{\Ga_2 \cdot
  \De_2} = X_{\Ga_1 \cdot \De_2} X_{\Ga_2 \cdot \De_1}$.  We call
these the {\em binomial relations}.  A special case are the {\em
  simple binomial relations} when $\vdeg \De_i = (1, 1, 1, 1, 0,
\dots, 0) = 1^4 0^{n-4}$, or some permutation thereof.  Examples 
are shown in Figures ~\ref{s5} and ~\ref{s835}.

(We have now defined all the relations relevant to the Main
Theorem~\ref{mainthm}, so the reader is encouraged to reread its
statement.)

In the even democratic case, the smallest quadratic binomial relations that are
not simple binomial relations appear for $n=12$, and
$$\vdeg \Ga_i = \vdeg \De_i = 
(1,1,1,1,1,1,0,0,0,0,0,0).$$
In the introduction, we asked if these quadratics are linear
combinations of the simple binomial relations. \label{sneaky1}\lremind{sneaky1}

By the Pl\"ucker relations, the binomial relations are generated by
those where the $\Ga_i$ and $\De_i$ are non-crossing, and similarly
for the simple binomial relations.  This restriction can be useful to
reduce the number of equations, but as always, symmetry-breaking
obscures other algebraic structures.  (Note that even though we may
restrict to the case where $\Ga_i$ and $\De_j$ are non-crossing, we
may not restrict to the case where $\Ga_i \cdot \De_j$ are
non-crossing, as the following examples with $n=5$ and $n=8$ show.)

As an example, consider the case $n=5$ (with the smallest democratic
linearization $(2,2,2,2,2)$).  One of the simple binomial
relations is shown in Figure~\ref{s5}. The building blocks $\Ga_i$ and
$\De_j$ are shown Figure~\ref{s5prep}.   In fact, these quadric relations
cut out $M_5$ in $\proj^5$, as can be checked directly, or as follows
from Theorem~\ref{mainthm}.  The $\SSS_5$-representation on the quadrics
is visible.

\begin{figure}
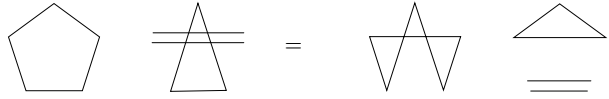

\begin{center}\include{s5}
\end{center}
\caption{A simple binomial relation for $n=5$ \label{s5}\lremind{s5}}
\end{figure}

\begin{figure}
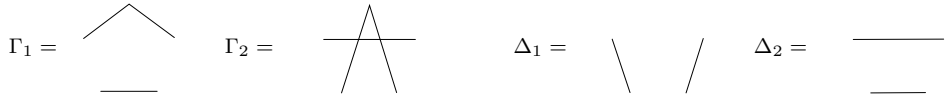

\begin{center}
\include{s5prep}
\end{center}
\caption{The building blocks of Figure~\ref{s5} \label{s5prep}\lremind{s5prep}}
\end{figure}

\point As a second example, consider $n=8$ (and the democratic
linearization $(1, \dots, 1)$).  Because there are $\binom 8 4 / 2 =
35$ ways of partitioning the $8$ vertices into two subsets of size
$4$, and each such partition gives one simple binomial relation (where
the $\Ga_i$ and $\De_j$ are non-crossing, see comments two paragraphs
previous), we have $35$ quadric relations on $M_8$, shown in
Figure~\ref{s835}.  \label{n8}\lremind{n8}

\begin{figure}
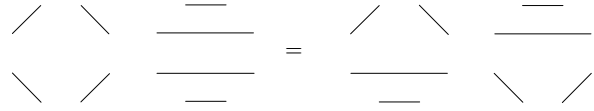

\begin{center}
\include{s835}
\end{center}
\caption{One of the simple binomial relations for $n=8$ points
\label{s835}\lremind{s835}}
\end{figure}

The space of quadric relations forms an irreducible $14$-dimensional
$\SSS_8$-representation, which we show by representation theory.
If $V$ is the vector space of quadratic relations, we have the
exact sequence
$$
0 \rightarrow V \rightarrow \Sym^2 H^0(M_8, \oh(1)) \rightarrow
H^0(M_8, \oh(2)) \rightarrow 0$$
of representations.  By counting
non-crossing graphs, we can calculate $h^0(M_8, \oh(2)) = 91$
(the $8$th Riordan number or Motzkin sum) and 
$h^0(M_8, \oh(1))=C_4=14$ (the $4$th Catalan number), from which $\dim V =
14$.  As the representation $H^0(M_8, \oh(1))$ is identified,
we can calculate the
representation $\Sym^2 H^0(M_8, \oh(1))$, and observe that the only
$14$-dimensional subrepresentations it contains are irreducible.
(Simpler still is to compute the character of the $196$-dimensional
representation $H^0(M_8, \oh(1))^{\otimes 2}$, and decompose it into
irreducible representations, using Maple for example, finding that it decomposes into
representations of dimension $1+14+14+20+35+56+56$;
$\Sym^2 H^0(M_8, \oh(1))$ is of course a subrepresentation of this.)

As our quadric relations are nontrivial, and form an
$\SSS_8$-representation, we have given generators of the quadric
relations.  Necessarily they span the same vector space of the $14$
relations given in  \cite[\S 8.3.4]{hmsv}.  Our relations have the
advantage that the $\SSS_8$-action is clear, but the major disadvantage
that it is not a priori clear that the vector space they span has
dimension $14$.  We suspect there is an $\SSS_8$-equivariant description
of the linear relations among the generators, but we have been unable to
find one.

In \cite[\S 8.3.4]{hmsv} it was shown that the ideal of relations of
$M_8$ was generated by these fourteen quadrics, and hence by our $35$
simple binomial quadrics.  We will use this as the base case of our
induction later.

\epoint{The Segre cubic relation} Other relations are also clear from
this graphical perspective.  For example, Figure~\ref{snegre} shows an
obvious relation for $M_6$, which is well-known to be a cubic
hypersurface (the Segre cubic hypersurface, see for example
\cite[p.~17]{do} or \cite[Example~11.6]{d}).  As this is a nontrivial
cubic relation (this can be verified by writing it in terms of a
non-crossing basis), it must be the Segre cubic relation.
Interestingly, although the relation is not $\SSS_6$-invariant, it
becomes so modulo the Pl\"ucker relations \eqref{plucker}.  Note that
there are no (nontrivial) binomial relations for $M_6$ (which is cut
out by this cubic), so the Segre relation cannot be in the ideal
generated by the binomial relations.

\epoint{Remark: Segre cubic relations for $n \geq 8$} There are
analogous cubic relations for $n \geq 8$, by simply adding other
vertices.  The $n=8$ case is given in Figure~\ref{snegreplus}.  For $n
\geq 8$, these Segre cubic relations lie in the ideal generated by the
simple binomial relations. We will use this in the proof of
Theorem~\ref{mainthm}.  This follows from the case $n=8$, which can be
verified in a couple of ways.  As stated above, \cite{hmsv} shows that
the ideal cutting out $M_8$ is generated by the fourteen quadrics of
\cite[\S 8.3.4]{hmsv}, which by \S \ref{n8} is the ideal generated by
the simple binomial relations, and the cubic lies in this ideal.  One
may also verify that the Segre relation lies in the ideal generated by
the fourteen quadrics by explicit calculation (omitted here).
\label{segresurprise}\lremind{segresurprise}

\begin{figure}
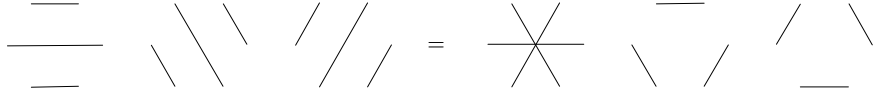

\begin{center}
\include{snegre}
\end{center}
\caption{The Segre cubic relation (graphical version)
\label{snegre}\lremind{snegre}}
\end{figure}

\begin{figure}
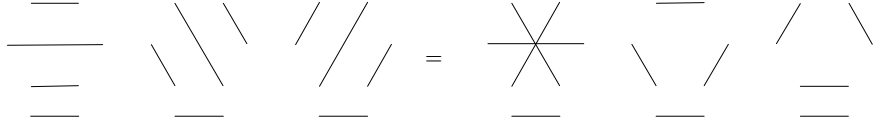

\begin{center}
\include{snegreplus}
\end{center}
\caption{The Segre cubic relation for $n=8$
\label{snegreplus}\lremind{snegreplus}}
\end{figure}

\epoint{Other relations} There are other relations, that we will not
discuss further.  
For example,  consider the democratic case for $n$ even.
Then $\SSS_n$ acts on the set of graphs.  Choose any graph $\Ga$.
Then 
$$
\sum_{\si \in \SSS_n} \operatorname{sgn}(\si) X_{\si(\Ga)}^i = 0
$$
is a  relation for $i$ odd and $1<i< n-1$.  Reason:
substituting for $X$'s in terms of $p$'s (or more correctly the $u$'s
and $v$'s) using \eqref{snack} to obtain an expression $E$, and
observing that $\SSS_n$ acts oddly on $E$, we see that we must obtain
a multiple of the Vandermonde, which has degree $(n-1, \dots, n-1) > \deg E$.
Hence $E=0$.  It is not clear that this is a nontrivial relation,
but it appears to be so in small cases.  In particular, the case $n=6$,
$i=3$ is the Segre cubic relation.
In the introduction, we asked if these relations
for $n=10$ lie in the ideal generated by the simple binomial
quadric relations. \label{sneaky2}\lremind{sneaky2}

\bpoint{Degree of the GIT quotient}
\label{degree}\lremind{degree}As an application of these coordinates, we compute the degree of all
$M_{\vw}$.  For example, we will use this to verify that the
degree is $1$ when $\left| \vw \right| = 6$ and $\vw \neq (1,\dots,
1)$, although this can also be done directly.

  We would like to intersect the moduli space
$M_{\vw}$ with $n-3$ coordinate hyperplanes of the form  $X_{\Ga}=0$ and count the
number of points, but these hyperplanes will essentially never
intersect properly.  Instead, we note that each hyperplane $X_{\Ga}=0$
is reducible, and consists of a finite number of components of the
form $M_{\vw'}$ where the number of points $\# \vw'$ is $n-1$.  We can compute the multiplicity with which each of these
components appears.  The algorithm is then complete, given  the base case $n=4$.
Here, more precisely, is the algorithm.  

(a) (trivial case) If $n=3$, the moduli space is a point, so the
degree is $1$.  

(b) (base case) If $\vw = (d,d,d,d)$, then $\deg M_{\vw} = d$, as the
moduli space is isomorphic to $\proj^1$, embedded by the $d$-uple
Veronese.  (This may be seen by direct calculation, or by noting that
a base-point-free subset of those variables of degree $(d,d,d,d)$ are
``$d$th powers'' of variables of degree $(1,1,1,1)$.)

(c) If $n>4$ and $\vw$ satisfies $w_j + w_k \leq \sum w_i/2$ for all
$j$, $k$, we choose any $\Ga$ of weight $\vw$.  We can understand the
components of $X_{\Ga}=0$ by considering the
morphism $\pi: (\proj^1)^n - U_{\vw} \rightarrow [ X_{\Ga} ]_{\vdeg
  \Ga = \vw}$, where $U_{\vw}$ is the unstable locus.  Directly from
the formula for $X_{\Ga}$, we see that for each pair of vertices $j$,
$k$ with an edge joining them, such that $w_j+w_k < \sum w_i / 2$,
there is a component that can be interpreted as $M_{\vw'}$, where
$\vw'$ is the same as $\vw$ except that $w_j$ and $w_k$ are removed,
and $w_j+w_k$ is added (call this $w_0$ for convenience).  We
interpret this as removing vertices $j$ and $k$, and replacing them
with vertex $0$.  This component corresponds to the divisor
\lremind{snuffle}
\begin{equation}
\label{snuffle}
(u_j v_k - u_k
v_j)^{m_{jk}}=0
\end{equation} on the source of $\pi$, where $m_{jk}$ is the number of
edges joining vertices $j$ and $k$.  If $\De$ is the reduced
version of this divisor, $u_j v_k - u_k v_j=0$, then the correspondence
between between $\De \rightarrow M_{\vw}$ and
$(\proj^1)^{n-1} - U_{\vw'} \rightarrow M_{\vw'}$ is as follows.
For each $\Ga'$ of degree $\vw'$, we lift $X_{\Ga'}$ to any $X_{\Ga}$ where
$\Ga$ is a graph on $\{ 1, \dots, n \}$ of degree $\vw$ whose ``image'' in $\{ 1, \dots, n \} \cup \{ 0 \}
\setminus \{ j,k \}$  is $\Ga'$.  (In other words, to $w_j$ of the
$w_0$
edges meeting vertex $0$ in $\Ga'$, we associate edges meeting vertex
$j$ in $\Ga$, and similarly with $j$ replaced by $k$.)
If $\Ga''$ is any other lift, then $X_{\Ga} = \pm  X_{\Ga''}$ on $\De$, 
because using the Pl\"ucker relations, $X_{\Ga} \pm X_{\Ga''}$  can be expressed
as a combination of variables containing edge $jk$, which all vanish on $\De$.

From \eqref{snuffle}, the multiplicity with which this component
appears is $m_{ij}$, the number of edges joining vertices $j$ and $k$.

If $w_j + w_k = \sum w_i/2$, then $M_{\vw'}$ is a strictly semistable
point, and of dimension $0$ smaller than $\dim M_{\vw} - 1$, and hence is
not a component.  (Our base case is $n=4$, not $3$, for this reason.)

(d) If $n \geq 4$ and there are $j$ and $k$ such that $w_j + w_k > \sum w_i/2$,
then the rational map
$(\proj^1)^4 \dashrightarrow M_{\vw}$ has a base locus.  Any
graph $X_{\Ga}$ of degree $\vw$ necessarily contains a copy of edge
$jk$, so $(u_j v_k - u_k v_j)$ is a factor of any of the $X_{\Ga}$.
Hence $M_{\vw}$ is naturally isomorphic to $M_{\vw - e_j - e_k}$,
so we replace $\vw$ by $\vw - e_j - e_k$, and repeat the process.
Note that if $n=4$, then the final resulting quadruple must be of the
form $(d,d,d,d)$.

For example, $\deg M_4 = 1$, $\deg M_6 = 3$, $\deg M_8 = 40$, and
$\deg M_{10}= 1225$ were computed by hand.  (This appears to be
sequence A012250 on Sloane's {\em On-line encyclopedia of integer
  sequences} \cite{sloane}.)  The calculations $\deg M_6=3$ and $\deg
M_{2,2,2,2,2} = 5$ are shown in Figure~\ref{m6}. and \ref{m5}
respectively.  At each stage, $\vw$ is shown, as well as the $\Ga$
used to calculate the next stage.  In these examples, there is
essentially only one such $\vw'$ at each stage, but in general there will be
many. The vertical arrows correspond to identifying components of
$X_{\Ga}$ (step (c)).  The first arrow in Figure~\ref{m6} is labeled
$\times 3$ to point out the reader that the next stage can be obtained
in three ways.  The degrees are obtained inductively from the bottom
up.  (The reader is encouraged to show that $\deg M_8 = 40$, and that this algorithm indeed gives $\deg M_{d \vw} = d^{n-3} \deg
M_{\vw}$.)

\begin{figure}
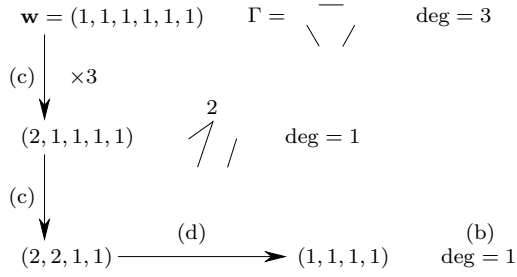

\begin{center}
\include{m6}
\end{center}
\caption{Computing $\deg M_6 =  3$ (recall that $M_6$
is the Segre cubic threefold $\Segre$)
\label{m6} \lremind{m6} }
\end{figure}

\begin{figure}
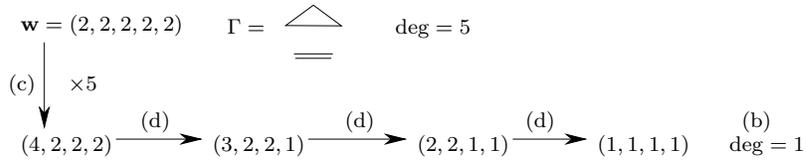

\begin{center}
\include{m5}
\end{center}
\caption{Computing $\deg M_{(2,2,2,2,2)} =  5$
using an inconvenient choice of $\Ga$
(recall that $M_{(2,2,2,2,2)}$ is a degree $5$ del Pezzo surface)
\label{m5} \lremind{m5} }
\end{figure}

\bpoint{Reduction of the Main Theorem~\ref{mainthm} to the equilateral case}
\label{undem}\lremind{undem}We next show that the Main Theorem~\ref{mainthm} in the equilateral
case (when $\vw = (1, \dots, 1)$) implies the Main
Theorem~\ref{mainthm} in general.  The argument is similar in spirit to
our proof of Kempe's Theorem ~\ref{kempe}. 

Consider the commutative diagram:
$$\xymatrix{
0 \ar[r] & I_n \ar[r] \ar[d]^{\al} & 
\oplus_{\Gamma} \Z X_\Gamma  \ar[r]^{\phi} \ar[d]^{\be} &  R_n \ar[d]^{\ga} \ar[r] & 0 \\
0 \ar[r] &
I_{\vw} \ar[r] & \oplus_\Omega\Z X_\Omega  \ar[r]^{\psi} & R_{\vw} \ar[r] & 0, \\
}$$
where $R_n$
(resp.\ $R_{\vw}$) 
is the ring of invariants for $1^n$ (resp.\ $\vw$),
the $\Gamma$'s range over matchings of $n$ vertices, $n = | \vw |$ is even, 
and the $\Omega$'s range over multi-valence $\vw$ graphs.  The map $\be$ 
takes $X_\Gamma$ to $X_\Omega$ where $\Omega$ is given by identifying vertices of $\Gamma$ within the same clump; if a loop is introduced then it maps to zero.  
(We take the clumps to be subsets of adjacent vertices in the $n$-gon.)

\tpoint{Theorem}  {\em The map $\al : I_n \to I_{\vw}$ is surjective.}

\bpf 
It is clear that $\be$ is surjective (and hence  $\ga$ too).
Thus by the five lemma, it suffices to prove that $\ker \be$ 
surjects onto $\ker \ga$.

We have that $R_n = \oplus_G \Z \cdot \phi(X_G)$, where $G$ ranges
over regular non-crossing graphs on $\{1,\ldots,n\}$.  (In other
words, such $\phi(X_G)$ form a $\Z$-basis of $R_n$.)  Similarly
$R_{\vw} = \oplus_H \Z \cdot \psi(X_H)$ as $H$ ranges over
non-crossing graphs of multi-valence a multiple of $\vw$.  For each
$H$, there is exactly one non-crossing $G$ for which $\ga( \phi(X_G) )
= \psi(X_H)$.  Thus $\ker(\ga) = \oplus_G \Z \cdot\phi(X_G)$ where the
sum is over those non-crossing $G$ which contain least one edge which
connects two vertices in a single clump.  Fix such a $G$, with an edge
$a \to b$ where $a$ and $b$ are in the same clump.  We will show that
$X_G \in \ker \be$.

Partition $\{1,\ldots,n\}$ into two equal sized subsets $A$ and $B$
(``positive'' and ``negative'') such that $a \in A$, and $b \in B$.
As in the proof of Kempe's theorem~\ref{kempe}, we can write $X_G =
\sum_i \pm  X_{G_i}$, where the $G_i$ are (possibly
crossing) graphs each containing the edge $a \to b$.  (The process
described in the proof of Theorem~\ref{kempe} involves trading a pair
of edges, one positive and one negative, for two neutral edges.  No
neutral edges such as $a \rightarrow b$ are affected by this process.)

By applying Hall's marriage theorem (repeatedly) to each $G_i$, we can write
$X_{G_i} = \prod_{j = 1}^k \phi(X_{\Gamma_{i,j}}),$
where the $\Gamma_{i,j}$ are matchings, and $k = \deg(X_G)$.  (At least)
one  $\Gamma_{i,j}$ contains
the edge $a \to b$, so, $\be(\prod_{j = 1}^k X_{\Gamma_{i,j}}) = 0$. 
Hence $\be(X_G) = \be( \sum \pm  X_{G_i}) = 0$ as desired.
\epf

\bpoint{Verification of the Main Theorem~\ref{mainthm} in small cases}
\label{smallcases}\lremind{smallcases}The cases $\left| \vw \right| = 2$ and $\left| \vw
\right| = 4$ are trivial.

If $\left| \vw \right| = 6$ and $\vw \neq
(1, \dots, 1)$, then $\vw = (3,2,1)$, $(2,2,2)$, $(2,2,1,1)$, or
$(2,1,1,1,1)$.  The first two cases are points, and the next two cases
were verified to have degree $1$ in \S \ref{degree} (see
Figure~\ref{m6}).

The case $\vw= (1,1,1,1,1,1,1,1)$ was verified in \S \ref{n8}, so by \S \ref{undem},
the case $\left| \vw \right|=8$ follows. 

Thus the cases $\left| \vw \right| \geq 10$ remain.

\section{An analysis of a neighborhood of a strictly semistable point}

\label{sss}\lremind{sss}We now show the result in a neighborhood of a strictly semistable
point, in the equilateral case $w=1^{n=2m}$, in characteristic $0$, by
explicitly describing an affine neighborhood of such a point.  This
affine neighborhood has a simple description: it is the space of
$(m-1) \times (m-1)$ matrices of rank at most $1$, where no entry is
$1$ (Lemma~\ref{rankone}).  The strictly semistable point corresponds
to the zero matrix.

\bpoint{The Gel'fand-MacPherson correspondence: the moduli space
  as a quotient of the Grassmannian}
We begin by recalling the Gel'fand-MacPherson correspondence, an
 alternate description of the
moduli space.
The Pl\"ucker embedding of the Grassmannian $G(2,n) \hookrightarrow
\proj^{\binom n 2 -1}$ is via the line bundle $\oh(1)$ that is the
positive generator for $\Pic G(2,n)$.  This generator may be
described explicitly as follows.  Over $G(2,n)$, we have a
tautological exact sequence of vector bundles \lremind{taut}
\begin{equation}
\xymatrix{ 0 \ar[r] & S  \ar[r] & \oh^{\oplus n} \ar[r] & Q  \ar[r] & 0} 
\label{taut}
\end{equation}
where $S$ is the tautological rank $2$ subbundle (over $[\Lambda ] \in
G(2,n)$, it corresponds to $\Lambda \subset \C^n$), and $Q$ is the
tautological rank $n-2$ quotient bundle.  Then $\wedge^2 S = \oh(1)$ is a line
bundle, and is the dual to $\oh(1)$.  Dualizing \eqref{taut} we 
get a map $\wedge^2 \oh^{\oplus n} \rightarrow \wedge^2 S^*$.  
Then $\wedge^2 S^*$  can be easily checked to be be generated
by these global sections.
We call these sections $s_{ij}$, and note that they satisfy the following
relations: the sign relations $s_{ij} = -s_{ji}$ (so $s_{ii} =0$)
inherited from $\wedge^2 \oh^{\oplus n}$, and the Pl\"ucker relations
$$
s_{ij} s_{kl} - s_{ik} s_{jl} - s_{jk} s_{il} = 0.$$ These equations
cut out the Grassmannian in $\proj^{\binom n 2 -1}$.

The connection to $n$ points on $\proj^1$ is as follows.  Given a
general point of the Grassmannian corresponding to the subspace
$\Lambda$ of $n$-space, we obtain $n$ points on $\proj^1$, by
considering the intersection of $\Lambda$ with the $n$ coordinate
hyperplanes, and projectivizing.  This breaks down if $\Lambda$ is
contained in a coordinate hyperplane.  \cut{(Aside: if $\Lambda$ is
  contained in coordinate hyperplane $H_i$, then $s_{ij}=0$ for all
  $j$, so these very bad points are discarded in the definition of
  $V_P$ below.)}  (The point $[\Lambda ]$ is GIT-stable if the
resulting $n$ points in $\proj^1$ are GIT-stable, and similarly for
semistable.  We recover the cross-ratio of four points via
$s_{ij}s_{kl}/s_{il}s_{jk}$.)

Let $D(s_{1n})$ be the distinguished open set where $s_{1n} \neq 0$.
(In the correspondence with marked points, this corresponds to the locus
where the first point is distinct from the last point.)
Then $D(s_{1n})$ is isomorphic to $\A^{2(n-2)}$, with good
coordinates as follows.  Given $\Lambda \notin D(s_{1n})$,
choose a basis for $\La$, written as a $2 \times n$ matrix.  As $\Lambda
\notin D(s_{1n})$, the first and last columns are linearly
independent, so up to left-multiplication by $GL(2)$, there is a
unique way to choose a basis where the first column is $\left[
\begin{array}{c} 0 \\ 1 \end{array} \right]$ and the last column is $
\left[ 
  \begin{array}{c} 1 \\ 0 \end{array}\right]$.  We choose the
``anti-identity'' matrix rather than the identity matrix, because we
will think of the first column as $[0;1] \in \proj^1$, and the last
column as $[1;0]$.  (Another interpretation is as follows.  If $\La$
is interpreted as a line in $\proj^{n-1}$, and $H_1$, \dots, $H_n$ are
the coordinate hyperplanes, then if $\La$ does not meet $H_1 \cap
H_n$, then it meets $H_1$ at one point of $H_1 - H_1 \cap H_n \cong
\C^{n-2}$ and $H_n$ at one point of $H_n - H_1 \cap H_n \cong
\C^{n-2}$, and $\Lambda$ is determined by these two points. The
coordinates on the first space are the $x$'s, and the coordinates on
the second are the $y$'s.)

Thus if the $2 \times n$ matrix is written
$$
\left[ \begin{array}{cccccc}
0 & x_2 & x_3 & \cdots & x_{n-1} & 1 \\
1 & y_2 & y_3 & \cdots & y_{n-1} & 0 
\end{array} \right]
$$
then we have coordinates $x_2$, \dots, $x_{n-1}$, $y_2$, \dots, $y_{n-1}$
on our affine chart.  For convenience, we define $x_1=0$, $y_1=1$, $x_n=1$,
$y_n=0$.  

Under the trivialization $(\oh(1), s_{1n})|_{D(s_{1n})} \cong (\oh,
1)|_{D(s_{1n})}$, in these coordinates, the section $s_{ij}$ may be
interpreted as $$ s_{ij} = x_j y_i - x_i y_j.
$$
We can use this to immediately verify the Pl\"ucker relations.
We also recover the $x_i$ and $y_j$ from the sections via\lremind{xiyj}
\begin{equation}\label{xiyj}
x_i = s_{1i} / s_{1n} \quad \quad
y_j = s_{jn} / s_{1n}.\end{equation}

The Grassmannian has dimension $2(n-2) = 2n-4$.  To obtain our moduli
space, we take the quotient of $G(2,n)$ by the maximal torus $T
\subset SL(2,n)$, which has dimension $n-1$.  (Thus as expected the
quotient has dimension $n-3$.)  We will write elements of this maximal
torus as $\lambda = (\lambda_1, ..., \lambda_n)$.  To describe the
linearization, we must describe how $\lambda$ acts on each $s_{ij}$:
$\lambda_i$ acts on $s_{ij}$ with weight $1$, and on the rest of the
$s_{ij}$'s by weight $0$.  This action certainly preserves our
relations.

Then we can see how to construct the quotient as a $\Proj$: the terms
that have weight $(d,d,\dots, d)$ correspond precisely to $d$-regular graphs on our $n$
vertices.  
Hence we conclude that this projective scheme is precisely
the GIT quotient of $n$ points on the projective line, as the graded
rings are the same.  This is the Gel'fand-MacPherson correspondence.
The relations we have described on our $X_{\Gamma}$ clearly
come from the relations on the Grassmannian.  (That is of course no
guarantee that we have them all!) 

\bpoint{A neighborhood of a strictly semistable point}
Let $\pi: G(2,2m)^{ss}
\rightarrow M_{\vw}$ be the quotient map.  Let $p$ be a strictly semistable 
point
of the moduli space $M_{\vw}$, without loss of generality the image of
$(0, \dots, 0, \infty, \dots, \infty)$.  We say an edge $ij$ on
vertices $\{ 1, \dots, 2m \}$ is {\em good} is $i \leq m < j$ (if it 
``doesn't connect two $0$'s or two $\infty$'s'').  We say
a graph on $\{ 1, \dots, 2m \}$ is {\em good} all of its edges are
good.  We say an edge or graph is {\em bad} if it is not good.  Let
$P$ be the set of good matchings of $\{ 1, \dots, 2m \}$.  Let
$$
U_P = \{ q \in M_{\vw} : X_{\Ga}(q) \neq 0 \text{ for all $\Ga \in P$} \}.$$
(In the dictionary to $n$ points on $\proj^1$, this corresponds
to the set where none of the first $m$ points is allowed to be the same
as any of the last $m$ points.)
Note that  $p \in P$, and $\pi^{-1} (U_P) \subset D(s_{1, 2m})$.

\tpoint{Lemma}  {\em $U_P$ is an affine variety, with coordinate ring
generated by $W_{ij}$  and $Z_{ij}$ ($1 <  i \leq m <  j < 2 m$)
with relations \lremind{atomic1-2}\lremind{llemma}\label{llemma}
\begin{equation}
\label{atomic1}
W_{ij} W_{kl} = W_{il} W_{kj}\end{equation} (i.e.\  the matrix
$[W_{ij}]$ has rank $1$) and 
\begin{equation}Z_{ij}(W_{ij}-1) = 1 \label{atomic2}
\end{equation}
(i.e.\ the matrix $[W_{ij}]$ has no entry $1$).}\label{rankone}\lremind{rankone}

This has a simple interpretation: $U_P$ is isomorphic to the space of
$(m-1) \times (m-1)$ matrices of rank at most $1$, where each entry
differs from $1$, and $p$ is the unique singular point, corresponding
to the zero matrix. 
(We remark that this is the cone over the Segre
embedding of $\proj^{m-1} \times \proj^{m-1}$.)
 Hence we have described a neighborhood of the
singular point rather explicitly.

\bpf
Let $V_P = \{ [\Lambda] \in G(2,2m) : s_{ij}(\Lambda) \neq 0
\text{ for all $i \leq m < j$ }\}$, so $V_P = \pi^{-1}(U_P)$.
Then $V_P$ is an open subset of $D(s_{1, 2m})$.
In terms of the coordinates on $D(s_{1, 2m}) \cong \C^{4m-4}$
described above,  $V_P$ is described by\lremind{VPeqns}
\begin{equation}
\label{VPeqns}
x_j y_i - x_i y_j \neq 0, \quad \quad\quad x_j \neq 0, \quad \quad\quad y_i \neq 0
\end{equation}
for $i \leq m < j$.
Let $T$ be the 
maximal torus $T \subset SL(2m)$.
By our preceding discussion, using \eqref{xiyj}, $\la_k$ ($1<k<2m$)
acts on $x_k$ and $y_k$ with weight $1$, and on the other
$x_l$ and $y_l$'s with weight $0$.
The torus $\la_1$ acts on $x_l$ with weight $0$, and on $y_l$
with weight $-1$.
The torus $\la_{2m}$ acts on $x_l$ with weight $-1$, and on $y_l$
with weight $0$.

We analyze the quotient $V_P / T$ by writing $T$ as a product, $T= T'
T''$, and restricting to invariants of $T'$, then of $T''$.  Let $T'$
be the subtorus of $T$ such that $\la_1 = \la_{2m}$, and $T''$ be the
subtorus $(\la, 1, \dots, 1, \la^{-1})$.  We take the quotient first
by $T'$.  It is clear that the invariants are given by $u_i = x_i/y_i$
for $1 <i \leq m$ and $v_j = y_j/x_j$ for $m+1 \leq j < 2m$.  From
\eqref{VPeqns}, the quotient of $V_P$ is cut out by the inequality $u_i v_j-1 \neq 0$.

$T''$ acts on this quotient as follows: $\la$ acts on $u_i$ by weight
$-2$, and $v_j$ by weight $2$.  The invariants of this quotient by $T''$
are therefore  generated by $u_i v_j$ and $1/ (u_i v_j-1)$.

Hence if we take $W_{ij} = u_i v_j$, the invariants are generated
by the $W_{ij}$, subject to the relations that the matrix $[W_{ij}]$
has rank $1$, and also $W_{ij} \neq 1$.
\epf

Now let $I \subset \C[  \{ X_{\Ga} \} ]$ be the ideal
of relations of the invariants of $M_{\vw}$, and let
$I_V$ be the ideal generated by the linear Pl\"ucker relations
and the simple binomial relations.  We have already
shown that $I_V \subset I$.

Let $\SP$ be the multiplicative system of monomials in $X_{\Ga}$
generated by those $X_{\Ga}$ where $\Ga \in P$.
\cut{remind{[$S_P$ is now $\SP$.]}}

\tpoint{Theorem} {\em If $n=2m \geq 8$, then $\SP^{-1} I_V = \SP^{-1}
  I$.  In other words, the sign, Pl\"ucker, and simple binomial
  relations cut out the moduli space on this open subset.}
\label{ssst}\lremind{ssst}

As the Main
Theorem~\ref{mainthm} is true for $n=2m=8$ (\S \ref{n8}), this
theorem holds in that ``base'' case.

\bpf
By $\Ga$ we will mean a general matching, and by $\De$, we will
mean a matching in $P$. 
We have a surjective
map 
$$\C[  \{ X_{\Ga} / X_{\De} \} ] / \SP^{-1} I_V \rightarrow 
\C[ \{ X_{\Ga} / X_{\De} \} ] / \SP^{-1} I$$ (that we wish to show is
an isomorphism), and Lemma~\ref{llemma} provides an isomorphism
$$\C[ \{ X_{\Ga} / X_{\De} \} ] / \SP^{-1} I \cong \oh(U_P) \cong 
\C[ \{ W_{ij} ,  Z_{ij} \} ]/ J_{WZ},$$
where $J_{WZ} \subset \C[ \{ W_{ij} \} , \{Z_{ij} \}]$ is
the ideal generated by the relations \eqref{atomic1} and \eqref{atomic2}.
\cut{[Ben's $X$'s are now $W$'s, and Ben's $Y$'s are now $-Z$'s.  This is
temporary notation, but I didn't want to have too many $X$'s.]}

By comparing the moduli maps, we see that this isomorphism is given
by\lremind{map}
\begin{equation}\label{map}
W_{ij} \mapsto \frac {X_{1i \cdot j(2m) \cdot \Ga}}
 {X_{1j \cdot i(2m) \cdot \Ga}}, \quad \quad
Z_{ij} \mapsto \frac  {X_{1j \cdot i(2m) \cdot \Ga}}
{X_{1(2m) \cdot ij \cdot \Ga}}
\end{equation}
where $\Ga$ is any matching on $\{ 1, \dots, 2m \} - \{1, i, j, 2m \}$
such that $1j \cdot i(2m) \cdot \Ga \in P$.  (By the simple binomial
relations, this is independent of $\Ga$.)  The description of the
isomorphism in the reverse direction is not so pleasant, and we will
spend much of the proof avoiding describing it explicitly.

We thus have a surjective morphism
$$
\psi: \C[ \{ X_{\Ga}/X_{\De} \} ] \rightarrow \C[ \{ W_{ij}, Z_{ij} \} ]/J_{WZ}
$$
\cut{[Ben's $\Psi$.]}  whose kernel is $\SP^{-1} I$, which contains $\SP^{-1}
I_V$.  We wish to show that the kernel is $\SP^{-1} I_V$.  We do this as
follows. 
For each $1 < i \leq m < j < 2m$, fix a matching
$\Ga_{i,j}$ on $\{ 1, \dots, 2m \} - \{ 1, i, j , 2m \}$ 
so that $1j \cdot i (2m) \cdot \Ga_{i,j} \in P$.
Consider the subring of $\C [ \{ X_{\Ga} / X_{\De} \}]$
generated by 
$$w_{ij} = \frac {X_{1i \cdot j(2m) \cdot \Ga_{i,j}}}
 {X_{1j \cdot i(2m) \cdot \Ga_{i,j}}}, \quad \quad
z_{ij} = \frac  {X_{1j \cdot i(2m) \cdot \Ga_{i,j}}}
{X_{1(2m) \cdot ij \cdot \Ga_{i,j}}}.
$$
(Compare this to \eqref{map}.)
Call this subring $\C[ \{ w_{ij}, z_{ij} \} ] / J_{wz}$.

The proof consists of two steps.  {\em Step 1.}  We show that {\em
  any} element of $\C[ \{ X_{\Ga} / X_{\De} \} ]$ differs from an
element of $\C[ \{ w_{ij}, z_{ij} \} ]/ J_{wz}$ by an element of
$\SP^{-1} I_V$.  We do this in several smaller steps.  {\em Step 1a.}
We show that any $X_{\Ga}/X_{\De}$ can be written as a linear
combination of $X_{\De'}/X_{\De}$ (where $\De'$ is also good).  {\em
  Step 1b.}  We show that any such $X_{\De'}/X_{\De}$ may be expressed
(modulo $\SP^{-1} I_V$) in terms of $X_{ ik \cdot jl \cdot \Ga} /
X_{il \cdot jk \cdot \Ga}$, where $i,j \leq m < k,l$, and $\Ga$ is
good.  {\em Step 1c.}  We show that any such expression can be written
(modulo $\SP^{-1} I_V$) in terms of $w_{ij}$ and $z_{ij}$, i.e.\
modulo $\SP^{-1} I_V$, any such expression lies in $\C [ \{ w_{ij},
z_{ij} \} ] / J_{wz}$.

{\em Step 2.}  The kernel of the map $\psi: \C [ \{ w_{ij}, z_{ij} \}
] / J_{wz} \rightarrow \C[ \{ W_{ij}, Z_{ij} \} ]/ J_{WZ}$ (given by
$w_{ij} \mapsto W_{ij}$, $z_{ij} \mapsto Z_{ij}$) lies in $\SP^{-1}
I_V$. 

We now execute this strategy.

{\em Step 1a.}  We first claim that 
$X_{\Ga} / X_{\De}$ 
($\De \in P$)
is a linear combination of units 
$X_{\De'} / X_{\De}$ 
(i.e.\ 
$\De' \in P$) modulo the Pl\"ucker relations (the linear relations, 
which are in
$\SP^{-1} I_V$). 
We prove the result by induction on the number of bad
edges.  The base case (if all edges of $\Ga$ are good, i.e.\ $\Ga \in
P$) is immediate.  Otherwise, $\Ga$ has at least two bad edges, say
$ij$ and $kl$, where $i,j \leq m < k,l$.  
Then $X_{\Ga} = \pm X_{\Ga -
  \{ ij, kl \} + \{ ik, jl \}} \pm X_{\Ga - \{ ij, kl \} + \{ il, jk
  \}}$ is a Pl\"ucker relation, and the latter two terms have two fewer
bad edges, completing the induction.

{\em Step 1b.}  We show that any element $X_{\De'}/X_{\De}$
of $\C[ \{ X_{\Ga} / X_{\De} \} ]$ ($\De'$ good)
is congruent (modulo $\SP^{-1} I_V)$ to an element of the form
$X_{ ik \cdot jl \cdot \Ga} / X_{il \cdot
  jk \cdot \Ga}$, where $i,j \leq m < k,l$, and $\Ga$ is good. 
We prove this by induction on $m$.  If $m=4$, the result
is true (\S \ref{n8}).  Assume now that $m>4$.
If $\De'$ and $\De$ share an edge $e$, then let $\overline{\De'}$
and $\overline{\De}$ be the graphs on $2m-2$ vertices
obtained by removing this edge $e$.  Then by the inductive
hypothesis, the result holds for $X_{\overline{\De'}}/ X_{\overline{\De}}$.
By taking the resulting expression, and ``adding edge $e$ to the subscript
of each term'', we get an expression for $X_{\De'} / X_{\De}$.
Finally, if $\De'$ and $\De$ share no edge, suppose in $\De'$,
$1$ is connected to $(m+1)$; in $\De$, $1$ is connected to $(m+2)$;
and in $\De'$, $(m+2)$ is connected to $2$.
This is true after suitable reordering.
Say $\De' = 1 (m+1) \cdot 2(m+2) \cdot \Ga'$ and $\De = 1(m+2) \cdot \Ga$.
Then
$$
\frac 
{X_{\De'}}
{X_{\De}}
=
\frac {X_{1(m+2) \cdot 2(m+1) \cdot \Ga'}} { X_{1(m+2) \cdot \Ga}} 
\cdot
\frac {X_{1(m+1) \cdot 2(m+2) \cdot \Ga'}} 
{X_{1(m+2) \cdot 2(m+1) \cdot \Ga'}} .
$$
For each factor of the right side, the numerator and the denominator
``share an edge'', so we are done.

{\em Step 1c.}  We next show that any such can be written (modulo
$\SP^{-1} I_V$) in terms of $w_{ij}$ and $z_{ij}$, i.e.\ modulo
$\SP^{-1} I_V$ lies in $\C [ \{ w_{ij}, z_{ij} \} ] / J_{wz}$.  If
$2m=8$, the result again holds (\S \ref{n8}).  Assume now that
$2m>8$.  Given any $X_{ik \cdot jl \cdot \Ga}/ X_{ij \cdot kl \cdot
  \Ga}$ as in Step 1b, we will express it (modulo $\SP^{-1} I_V)$ in
terms of $w_{ij}$ and $z_{ij}$.  By the simple binomial relation
(i.e.\ modulo $\SP^{-1} I_V$), we may assume that $\Ga$ is any good
matching on $\{ 1, \dots, 2m \} - \{ i,j,k,l\}$, and in particular
that there are edges $ab, cd \in \Ga$ such that $\{ 1, 2m \} \subset \{
a,b,c,d,i,j,k,l \}$.  Then the result for $m=4$ implies that we can
write $X_{ik \cdot jl \cdot ab \cdot cd}/ X_{ij \cdot kl \cdot ab
  \cdot cd}$ can be written in terms of $w_{ij}$ and $z_{ij}$ (in terms
of the ``$m=4$ variables'').  By taking this expression, and ``adding
in the remaining edges of $\Ga$'', we get the desired result for our
case.

{\em Step 2.}  We will show that the kernel of the map $\psi: \C [ \{
w_{ij}, z_{ij} \} ] / J_{wz} \rightarrow \C[ \{ W_{ij}, Z_{ij} \} ]/
J_{WZ}$ (given by $w_{ij} \mapsto W_{ij}$, $z_{ij} \mapsto Z_{ij}$)
lies in $\SP^{-1} I_V$.

In order to do this, we need only verify that the
relations \eqref{atomic1} and \eqref{atomic2}
are consequences of the relations in $\SP^{-1} I_V$.

We first verify \eqref{atomic1}.  By the simple binomial
relation, we may write\lremind{a3,a4}
\begin{equation}
W_{ij}=  \frac{ X_{1i \cdot j (2m) \cdot kl \cdot \Ga}}
{ X_{1j \cdot i (2m) \cdot kl \cdot \Ga}},
\quad \quad 
W_{kl} = \frac{ X_{1 k \cdot l (2m) \cdot ij \cdot \Ga}}
{ X_{1 l \cdot k (2m) \cdot ij \cdot \Ga}},
\label{a3}
\end{equation}
\begin{equation}
W_{il}=  \frac{ X_{1i \cdot l (2m) \cdot jk \cdot \Ga}}
{ X_{1l \cdot i (2m) \cdot jk \cdot \Ga}},
\quad \quad 
W_{jk} = \frac{ X_{1 j \cdot k (2m) \cdot il \cdot \Ga}}
{ X_{1 k \cdot j (2m) \cdot il \cdot \Ga}}.
\label{a4}
\end{equation}
We thus wish to show that modulo $\SP^{-1} I_V$, the product of the
terms in \eqref{a3} equals the product of the terms in \eqref{a4}.
Choose any edge $e \in \Ga$.  The analogous question with $m=4$, with
$\Ga-e$ ``removed from the subscripts'', is true (\S \ref{n8}).
Hence, by ``adding $\Ga-e$ back in to the subscripts'', we get
the analogous result here.

We next verify \eqref{atomic2}:
$$
1 - W_{i,j} = \frac{ X_{1j \cdot i (2m) \cdot \Ga_{i,j}}}
{ X_{1j \cdot i (2m) \cdot \Ga_{i,j}}}
-
\frac{ X_{1i \cdot j (2m) \cdot \Ga_{i,j}}}
{ X_{1j \cdot i (2m) \cdot \Ga_{i,j}}}
\equiv 
\frac{ X_{1 (2m)  \cdot  i j \cdot \Ga_{i,j}}}
{ X_{1j \cdot i (2m) \cdot \Ga_{i,j}}}
= 1/ Z_{i,j} \pmod {\SP^{-1} I_V}$$
where the  equivalence uses a linear Pl\"ucker relation.
\epf

\section{Proof of Main Theorem}
\label{proof}\lremind{proof}
We have reduced to the equilateral case $\vw = 1^n$, $n=2m$, where $n
\geq 10$.  

The reader will notice that we will use the simple binomial relations
very little.  In fact we just use the inductive structure of the
moduli space: given a matching $\De$ on $n-k$ of $n$ vertices ($4 \leq
k < n$), and a point $[ X_{\Ga} ]_{\Ga}$ of $V_n$, then either these
$X_{\Ga}$ with $\De \subset \Ga$ are all zero, or $[ X_{\Ga} ]_{\De
  \subset \Ga}$ satisfies the Pl\"ucker and simple binomial relations
for $k$, and hence is a point of $V_k$ if $k \neq 6$.  (The reader
should think of this rational map $[ X_{\Ga} ] \dashrightarrow
[X_{\Ga}]_{\De \subset \Ga}$ as a forgetful map, remembering only the
moduli of the $k$ points.)  In fact, even if $k=6$ (and $n \geq 8$),
the point must lie in $M_6$, as the simple binomial relations for
$n>6$ induce the Segre cubic relation (\S \ref{segresurprise}). The
central idea of our proof is, ironically, to use the case $n=6$, where
the Main Theorem~\ref{mainthm} doesn't apply.

We will call such $\De$, where the $X_{\Ga}$ with $\De \subset \Ga$
are not all zero and the corresponding point of $M_6$ is stable, a
{\em stable $(n-6)$-matching}.  One motivation for this definition is
that given a stable configuration of $n$ points on $\proj^1$, there
always exists a stable $(n-6)$-matching.  (Hint: Construct $\De$
inductively as follows.  We say two of the $n$ points are in the same
{\em clump} if they have the same image on $\proj^1$. Choose any $y$
in the largest clump, and any $z$ in the second-largest clump; $yz$ is
our first edge of $\De$.  Then repeat this with the remaining
vertices, stopping when there are six vertices left.)  Caution: This
is false with $6$ replaced by $4$!

Main Theorem~\ref{mainthm} is a consequence of the following two
statements, and Theorem~\ref{ssst}. Indeed, {\bf (I)} and {\bf (II)}
show Theorem~\ref{mainthm} set-theoretically, and scheme-theoretically
away from the strictly semistable points (in characteristic $0$), and
Theorem~\ref{ssst} deals with (a neighborhood of) the strictly
semistable points.

{\bf (I)} There is a natural bijection between points of $V_n$
with no stable $(n-6)$-matching, and strictly semistable points of $M_n$.

{\bf (II)} If $B$ is any scheme, there is a bijection between morphisms $B
\rightarrow V_N$ missing the ``no stable $(n-6)$-matching'' locus
(i.e.\ missing the strictly semistable points of $M_n$, by (I)) and
stable families of $n$ points $B \times \{ 1, \dots, n \} \rightarrow \proj^1$.
(In other words, we are exhibiting an isomorphism of functors.)

One direction of the bijection of (I) is immediate.
The next result shows the other direction.

\tpoint{Claim} {\em If $[ X_{\Ga} ]_{\Ga}$ is a
  point of $V_n$ ($n \geq 10$) having no stable $(n-6)$-matching, then
 $[ X_{\Ga} ]_{\Ga}$ is a strictly semistable stable point of
  $M_n$.}
\label{proofofi}\lremind{proofofi}

Several of the steps will be used in the proof of (II).  We give them
names so they can be referred to later.

\bpf Our goal is to produce a partition of $n$ into two subsets of size
$n/2$, such that the point of $M_n$ given by this partition is our
point of $V_n$.  Throughout this proof, partitions will be assumed to
mean into two equal-sized subsets.

We work by induction.  We will use the fact that the result is also
true for $n=6$ (tautologically) and $n=8$ (as $V_8=M_8$, \S \ref{n8}).

Fix a matching $\De$ such that $X_{\De} \neq 0$.  By
the inductive hypothesis,  each edge $xy$ yields a strictly
semistable point of $M_{n-2}$, and hence a partition of $\{1, \dots, n \}
- \{x,y \}$, by considering all matchings containing $xy$.  Thus for
each $xy \in \De$, we get a partition of $\{ 1, \dots, n\}- \{x,y \}$.
If $wx$, $yz$ are two edges of $\De$, then we get the same induced
partition of $\{ 1, \dots, n \} - \{ w,x,y,z \}$ (from the inductive
hypothesis for $n-4$), so all of these partitions arise from 
a single partition $\{ 1, \dots, n \} = S_0 \coprod S_1$.  

\epoint{$\De$ two-overlap argument} As this partition is determined
using any two edges of $\De$, we would get the same partition if we
began with any $\De'$ sharing two edges with $\De$, such that
$X_{\De'} \neq 0$.
\label{2overlapi}\lremind{2overlapi}

\noindent {\em Defining the map to $\proj^1$.} Define $\phi: S_0
\coprod S_1 = \{ 1, \dots, n \} \rightarrow \proj^1$ by $S_0
\rightarrow 0$ and $S_1 \rightarrow 1$.  For each matching $\Ga$,
define $X'_{\Ga}$ using these points of $\proj^1$.  Rescale (or
normalize) all the $X_{\Ga}$ so $X'_\De = X_{\De}$.  We will show that
$X'_{\Ga} = X_{\Ga}$ for all $\Ga$, which will prove
Claim~\ref{proofofi}.

\epoint{One-overlap argument}\label{1overlapi}\lremind{1overlapi} For
any $\Ga$ sharing an edge $xy$ with $\De$, $X'_{\Ga} = X_{\Ga}$, for
the following reason: $[ X_{\Xi} ]_{xy \in \Xi}$ lies in $M_{n-2}$ by
the inductive hypothesis, and this point of $M_{n-2}$ corresponds to
the map $\phi$ (as the partition $S_0 \coprod S_1$ was determined
using this point of $M_{n-2}$), so $[ X_{\Xi} ]_{xy \in \Xi} = [
X'_{\Xi} ]_{xy \in \Xi}$, and the normalization $X'_{\De} = X_{\De}
\neq 0$ ensures that $X'_{\Xi}= X_{\Xi}$ for all $\Xi$ containing
$xy$.

\epoint{Reduction to $\Ga$ with $X'_{\Ga} \neq 0$}
It suffices to prove the result for those graphs $\Ga$, all
of whose edges connect $S_0$ and $S_1$ (i.e.\ no edge is contained in
$S_0$ or $S_1$; equivalently, $X'_{\Ga} \neq 0$).  We show this by
showing that {\em any} $X_{\Ga}$ is a linear combination of such
 graphs, by induction on the number $i$ of edges of $\Ga$
contained in $S_0$ (= the number contained in $S_1$).  The base case
$i=0$ is tautological.  For the inductive step, choose an edge $wx \in
\Ga$ contained in $S_0$ and an edge $yz$ contained in $S_1$.  Then  the
Pl\"ucker relation  using  $\Ga$ and $wxyz$ 
(with appropriate signs depending on the directions of edges)
is
$$\pm X_{\Ga} \pm X_{\Ga - wx - yz + wy + xz} \pm X_{\Ga - wx-yz + wz
  + xy} = 0,$$
and both $\Ga - wx - yz + wy + xz$ and $\Ga - wx-yz +
wz + xy$ have $i-1$ edges contained in $S_0$, and the result follows.
\label{goodgraphs}\lremind{goodgraphs}

\epoint{$pqrs$ argument, first version} Finally, assume that $X'_{\Ga}
\neq 0$ and that $\Ga$ shares no edge
with $\De$.  See Figure~\ref{snape}. Let $qr$ be an edge of $\Ga$ (so $\phi(q) \neq \phi(r)$),
and let $pq$ and $rs$ be edges of $\De$ containing $q$ and $r$
respectively (so $\phi(p) \neq \phi(q)$ and
$\phi(r) \neq \phi(s)$).  Then $\phi(p) \neq \phi(s)$, as $\phi$ takes
on only two values.  Let $\De' = \De - pq - rs + qr + ps$, so
$X'_{\De'} \neq 0$ as $\phi(q) \neq \phi(r)$ and $\phi(p) \neq
\phi(s)$.  Then $X_{\De'} = X'_{\De'}$ by the one-overlap argument
\ref{1overlapi}, as $\De'$ shares an edge with $\De$ (indeed
all but two edges), so $X_{\De'} \neq 0$.  Hence by the $\De$ two-overlap
argument \ref{2overlapi}, $\De'$ defines the same partition $S_0
\coprod S_1$, and hence the same map $\phi: \{ 1, \dots, n \}
\rightarrow \proj^1$.  Finally, $\Ga$ shares an edge with $\De'$, so
$X'_{\Ga} = X_{\Ga}$ by the one-overlap argument \ref{1overlapi}.
\label{pqrsi}\lremind{pqrsi}

\cut{\epoint{$pqrs$ argument} Finally, assume that $X_{\Ga} \neq 0$
  (i.e.\ $\Ga$ is good in the sense of the previous paragraph), and
  that $\Ga$ shares no edge of $\De$.  Let $qr$ be an edge of $\Ga$,
  and let $pq$ and $rs$ be edges of $\De$ containing $q$ and $r$
  respectively (see Figure~\ref{snape}).  (i) If $\phi(p) \neq
  \phi(s)$, then let $\De' = \De - pq - rs + qr + ps$; then $\De'$
  defines the same partition $S_0 \coprod S_1$ by \S \ref{2overlapi}
  (the $\De$ two-overlap argument --- here we use that $\Ga$ is good,
  so that $X_{\De'} \neq 0$ remind{Make sure this is true!!}), and
  $\Ga$ and $\De'$ share an edge, so we are done by \S \ref{1overlapi}
  (the $\Ga/\De$ one-overlap argument).  (ii) If $\phi(p) = \phi(s)$,
  then $\phi(p) \neq \phi(r)$ (as $rs$ is an edge of $\De$).  Let $st$
  be the edge of $\Ga$ containing $s$.  (It is possible that $t=p$.)
  Let $\Ga' = \Ga - qr - st + rs + qt$ and $\Ga''= \Ga - qr - st + qs
  + rt$ be the other two terms in the Pl\"ucker relation for $\Ga$ for
  $pqrs$.  Then $\Ga'$ shares an edge with $\De$, so $X'_{\Ga'} =
  X_{\Ga}$ by \S \ref{1overlapi} (the $\Ga/\De$ one-overlap argument),
  and by applying (i) to $\Ga''$ (swapping the names of $r$ and $s$),
  $X'_{\Ga''} = X_{\Ga''}$, so by the Pl\"ucker relation, $X'_{\Ga} =
  X_{\Ga}$ as desired.  }

We have thus completed the proof of Claim
\ref{proofofi}.
\epf

\begin{figure}
\begin{center}
\include{snape}
\end{center}
\caption{The $pqrs$ argument (vertex $t$ is used in \S \ref{nicematch})
\label{snape}\lremind{snape}}
\end{figure}

{\bf Proof of (II).}  The result boils down to the following
desideratum: Given any $(n-6)$-matching $\De$ on some $\{ 1, \dots, n
\} - \{ a, b,c,d,e, f \}$, there should be a bijection between:
\begin{enumerate}
\item[(a)] morphisms $\pi: B \rightarrow V_n$ contained in the open subset
  where $\De$ is a stable $(n-6)$-matching, and
\item[(b)] stable families of points $\phi: B \times \{ 1, \dots, n \}
  \rightarrow \proj^1$ where $\phi|_{B \times \{ a, \dots,  f \}}$
is also a stable family, and for any edge $xy$ of $\De$,
$\phi|_{B \times \{x \}}$ does not intersect
$\phi|_{B \times \{ y \}}$.
\end{enumerate}

We have already described the map (b) $\Rightarrow$ (a) in
\S \ref{snooker}. We now describe the map (a) $\Rightarrow$ (b), and
verify that (a) $\Rightarrow$ (b) $\Rightarrow$ (a) is the identity.
(It will then be clear that (b) $\Rightarrow$ (a) $\Rightarrow$ (b) is
the identity: given a stable family of points parameterized by $B$, we
get a map from $B$ to an open subset of $M_n$, which is a fine moduli
space, hence (b) $\Rightarrow$ (a) is an injection.  The result then
follows from the fact that (a) $\Rightarrow$ (b) $\Rightarrow$ (a) is
the identity.)

We work by induction on $n$.  The case $n=8$ was checked earlier 
 (\S \ref{n8}).

\noindent 
{\em The map to $\proj^1$.} Given an element of (a), define a family
of $n$ points of $\proj^1$ (an element of (b)) as follows.  (i) $\phi:
B \times \{ a, \dots, f \} \rightarrow \proj^1$ is given by the
corresponding map $B \rightarrow M_6$.  (ii) If $yz$ is an edge of
$\De$, we define $B \times ( \{ 1, \dots, n \} - \{ y, z \})
\rightarrow \proj^1$ extending (i) by considering the matchings
containing $yz$, which by the inductive hypothesis give a point of
$M_{n-2}$.  (iii) The morphisms of (ii) agree ``on the overlap'', as
given two edges $wx$ and $yz$ of $\De$, we get $B \times ( \{ 1,
\dots, n \} - \{ w,x,y,z \} ) \rightarrow \proj^1$ by considering the
matchings containing $wx \cdot yz$, which by the inductive hypothesis
give a map to $M_{n-4}$.  Here we are using that $n \geq 10$; and if
$n=10$, we need the fact that the Segre cubic relation cutting out
$M_6$ is induced by the quadrics cutting out $M_n$ for $n \geq 8$
(Remark ~\ref{segresurprise}).  Thus we get a well-defined morphism
$\phi: B \times \{ 1, \dots, n \} \rightarrow \proj^1$.

\epoint{$\De$ two-overlap argument, cf.\ \S \ref{2overlapi}} If $\De'$
is another matching on $\{1, \dots, n \} - \{ a, \dots, f \}$ sharing
at least $2$ edges with $\De$, with $X_{\De' \cdot \Xi} \neq 0$ for
some matching $\Xi$ of $\{a, \dots, f \}$, we obtain the same $\phi$,
as $\phi$ can be recovered by considering only two edges of $\De$ when
using (ii).  \label{2overlap}\lremind{2overlap}

\noindent {\em Defining $X'$.}  Define $X'_{\Ga}$
for all matchings $\Ga$ using $\phi$ and the moduli morphism of eqn.\ 
\eqref{snack}.  The coordinates $X_{\Ga}$ are projective (i.e.\ the
set of $X_{\Ga}$ is defined only up to scalars); scale them so that
$X_{\De \cdot \Xi} = X'_{\De \cdot \Xi}$ for all matchings $\Xi$ of $\{
a, \dots, f \}$.  Note that if $xy$ is an edge of $\De$, then $\phi(x)
\neq \phi(y)$, as there exists a matching $\Xi$ of $\{ a, \dots, f \}$
such that $X'_{\De \cdot \Xi} \neq 0$.

The following result will confirm that (a) $\Rightarrow$ (b)
$\Rightarrow$ (a) is the identity, concluding the proof of (II).

\tpoint{Claim}  We have the equality $X_{\Ga} = X'_{\Ga}$ for all
$\Ga$. 

\bpf
This proof will occupy us until the end of \S \ref{final}.

\epoint{One-overlap argument} As in \S \ref{1overlapi}, the result
holds for those $\Ga$ sharing an edge $yz$ with $\De$: by considering
only those variables $X_{\Ga'}$ containing the edge $yz$ (including
both $X_{\Ga}$ and $X_{\De}$), we obtain a point of $M_{n-2}$.
This point of $M_{n-2}$ is the one given by $\phi$ (this was part of
how $\phi$ was defined), so $[X_{\Ga'}]_{yz \in \Ga'} =
[X'_{\Ga'}]_{yz \in \Ga'}$.  By choosing a matching $\Xi$ on $\{
a, \dots, f \}$ so that $X_{\De \cdot \Xi} \neq 0$, we have that $X_{\Ga} X'_{\De
  \cdot \Xi} = X_{\De \cdot \Xi} X'_{\Ga}$.  Using $X_{\De \cdot \Xi} =
X'_{\De \cdot \Xi} \neq 0$, we have $X_{\Ga} = X'_{\Ga}$, as desired.
\label{1overlap}\lremind{1overlap}

We now deal with the remaining case, where $\Ga$ and $\De$ share no
edge.

\epoint{Reduction to $\Ga$ with $X'_{\Ga} \neq 0$ (cf.\ \S
  \ref{goodgraphs})} It suffices to prove the result for those graphs
such that $X'_{\Ga} \neq 0$, or equivalently that for each edge
$xy$ of $\Ga$, $\phi(x) \neq \phi(y)$.  We show this by showing that
any $X_{\Ga}$ is a linear combination of such graphs, by induction on
the number of edges $xy$ of $\Ga$ with $\phi(x) = \phi(y)$.
For the purposes of this paragraph, call these {\em bad edges}.
The base
case $i=0$ is tautological.  For the inductive step, choose a bad edge
$wx \in \Ga$ (with $\phi(w) = \phi(x)$), and another edge $yz$ such that
$\phi(y), \phi(z) \neq \phi(w)$.  (Such an edge exists, as by
stability, less than $n/2$ elements of $\{1, \dots, n \}$ take the
same value in $\proj^1$.)  
Then  the
Pl\"ucker relation  using  $\Ga$ with respect to  $wxyz$ 
is
$$\pm X_{\Ga} \pm X_{\Ga - wx - yz + wy + xz} \pm X_{\Ga - wx-yz + wz
  + xy} = 0,$$
and both $\Ga - wx - yz + wy + xz$ and $\Ga - wx-yz +
wz + xy$ have at most $i-1$ bad edges, and the result follows.
\label{goodgraphs2}\lremind{goodgraphs2}

Recall that we are proceeding by induction.  We first deal with the
case $n \geq 14$, assuming the cases $n=10$ and $n=12$.  We will then
deal with these two stray cases.  This logically backward, but the $n
\geq 14$ case is cleaner, and the two other cases are similar but more
ad hoc.

\bpoint{\bf The case $n\geq 14$} {\em $pqrs$ argument, second version.} As $n \geq
14$, there is an edge $qr$ of $\Ga$ not meeting $abcdef$.  
See
Figure~\ref{snape}.
By \S
\ref{goodgraphs2}, we may assume $\phi(q) \neq \phi(r)$.  Let $pq$ and
$rs$ be the edges of $\De$ meeting $q$ and $r$ respectively (so $\phi(p)\neq \phi(q)$ and $\phi(r)\neq
\phi(s)$. 
 (i) (cf.\ the similar argument of \S \ref{pqrsi}) If
$\phi(p) \neq \phi(s)$, then let $\De' = \De - pq - rs + qr + ps$;
then $\De'$ defines the same family of $n$ points as $\De$ by the
two-overlap argument \S \ref{2overlap}, and $\Ga$ and $\De'$ share an
edge, so we are done by the one-overlap argument \S \ref{1overlap}.
(More precisely, this argument applies on the open subset of $B$ where
$\phi(p) \neq \phi(s)$.)  (ii) If $\phi(p) = \phi(s)$, then $\phi(p)
\neq \phi(r)$.  (More precisely, this argument applies on the open set
where $\phi(p) \neq \phi(r)$.) Let $st$ be the edge of $\Ga$
containing $s$.  (It is possible that $t=p$.)  Let $\Ga' = \Ga - qr -
st + rs + qt$ and $\Ga''= \Ga - qr - st + qs + rt$ be the other two
terms in the Pl\"ucker relation for $\Ga$ for $pqrs$.  Then $\Ga'$
shares edge $rs$ with $\De$, so $X'_{\Ga'} = X_{\Ga'}$ by the 
one-overlap argument \S \ref{1overlap}, and by applying (i) to $\Ga''$
(swapping the names of $r$ and $s$), $X'_{\Ga''} = X_{\Ga''}$, so by
the Pl\"ucker relation, $X'_{\Ga} = X_{\Ga}$ as
desired.\label{nicematch}\lremind{nicematch}

\bpoint{The cases $n=10$ and $n=12$}  We are assuming that $\Ga$ and
$\De$ share no edges. 
If there is
an edge of $\Ga$ not meeting $\{ a, \dots, f \}$, the $pqrs$-argument
\S \ref{nicematch} applies, so assume otherwise.  
 Divide $\{ 1, \dots, n \}$ into two subsets
$abcdef$ and $ghij$ (resp.\ $ghijkl$) if $n=10$ (resp.\ $n=12$), where
the edges of $\De$ are $gh$, $ij$, and (if $n=12$) $kl$.  By renaming
$abcdef$, we may assume the edges of $\Ga$ are $ag$, $bh$, $ci$, $dj$,
and either $ef$ (if $n=10$, see Figure~\ref{10case}) or $ek$ and $fl$
(if $n=12$, see Figure~\ref{12case}).

\begin{figure}
\begin{center}
\include{10case}
\end{center}
\caption{The problematic graphs for $n=10$
\label{10case}\lremind{10case}}
\end{figure}

\begin{figure}
\begin{center}
\include{12case}
\end{center}
\caption{The problematic graphs for $n=12$
\label{12case}\lremind{12case}}
\end{figure}

\point Suppose that $\phi(a) \neq \phi(b)$.  Note that we will only
use that $ag, bh \in \Ga$, $gh \in \De$, and $\phi(a) \neq \phi(b)$
--- we will use this argument again below.  There is a matching $\Xi$
of $cdef$ so that if $xy \in \Xi$, then $\phi(x) \neq \phi(y)$.  (This
is a statement about stable configurations of $6$ points on $\proj^1$:
if we have a stable set of $6$ points on $\proj^1$, then no three of
them are the same point.  Hence for any four of them $cdef$, we can
find a matching of this sort.)  Let $\De' = \Xi \cdot ab \cdot \De$.
Then by the simple binomial relations (our first invocation!)
$X_{\De'} X_{\Ga} = X_{\De' -ab- gh+ag+bh} X_{\Ga +ab+ gh-ag-bh}$ and
$X'_{\De'} X'_{\Ga} = X'_{\De' -ab- gh+ag+bh} X'_{\Ga +ab+ gh-ag-bh}$.
However, by the one-overlap argument \S \ref{1overlap}, $X_{\De'} =
X'_{\De'} \neq 0$ ($\De'$ and $\De$ share edge $ij$), $X_{\De' -ab-
  gh+ag+bh} = X'_{\De' -ab- gh+ag+bh}$ ($\De' - ab - gh + ag+bh$ and
$\De$ share edge $ij$), and $X_{\Ga +ab+ gh-ag-bh} = X'_{\Ga +ab+
  gh-ag-bh}$ ($\Ga + ab + gh -ag-bh$ and $\De$ share edge $gh$), so we
are done.
\label{pianotpib}\lremind{pianotpib}

We are left with the case $\phi(a) = \phi(b)$.

\point \label{n10} 
Suppose now that $n=10$.  As $\phi(a) = \phi(b)$, $\phi(b)$ is
distinct from $\phi(e)$ and $\phi(f)$ (as $\phi(a)$, \dots, $\phi(f)$ are
a stable set of six points on $\proj^1$).  By the Pl\"ucker relations
for $\Ga$ (using $agef$),
$$\pm X_{\Ga} \pm X_{\Ga - ag - ef + ae + gf} \pm X_{\Ga - ag - ef + af + eg} 
 = 0,$$
and similarly for the $X'$ variables.  By applying the argument
of \S \ref{pianotpib} with $e$ and $a$ swapped, we have
$X'_{\Ga - ag - ef + af + eg}  = 
X_{\Ga - ag - ef + af + eg}$,
and by applying the argument of \S \ref{pianotpib}
with $f$ and $a$ swapped, we have
$X'_{\Ga - ag - ef + ae + gf}= 
X_{\Ga - ag - ef + ae + gf}$, 
from which $X'_{\Ga} = X_{\Ga}$, concluding the $n=10$
case. \lremind{n10}

\point Suppose finally that $n=12$.  If $\phi(c) \neq \phi(d)$, we are done
(by the same argument as \S \ref{pianotpib}, with $ab$ replaced by
$cd$), and similarly if  $\phi(e) \neq \phi(f)$.  Hence the only
case left is if $\phi(a)=\phi(b)$, $\phi(c) = \phi(d)$, and
$\phi(e)=\phi(f)$, and (by stability of the 6 points $\phi(a),
\dots, \phi(f)$) these are three distinct points of $\proj^1$.
Consider the Pl\"ucker relation for $\Ga$ with respect to $bhci$.  One
of the other two terms is $\Ga -bh - ci + bi + ch$, and
$X'_{\Ga-bh-ci+bi+ch} = X_{\Ga-bh-ci+bi+ch}$ (by the same argument as
in \S \ref{pianotpib}, as $\phi(a) \neq \phi(c)$).  We thus have to
prove that 
$X_{\Ga'} = X'_{\Ga'}$ for the third term in the Pl\"ucker relation, where
$$\Ga' = ag \cdot bc \cdot hi \cdot dj
\cdot ek \cdot fl.$$
For this, apply the argument of \S \ref{n10}, with $abghef$ replaced by $felkbc$ respectively.  \epf
\label{final}\lremind{final}

\cut{Old end (as of Aug 31 I think):  $dejk$.  One of the
other terms is $ag \cdot bc \cdot hi \cdot dk \cdot ej \cdot fl$, and
$X'_{ag \cdot bc \cdot hi \cdot dk \cdot ej \cdot fl} = X_{ag \cdot bc
  \cdot hi \cdot dk \cdot ej \cdot fl}$ (by the same argument as in
\S \ref{pianotpib}), so we must prove the result for the third
term in the Pl\"ucker relation:
$$\Ga'' = ag \cdot bc \cdot hi \cdot de \cdot jk \cdot fl$$
(see Figure~\ref{final12}).
Next
apply the Pl\"ucker relation for $\Ga''$ using $hide$.  Each of the
other two terms will have $h$ connected to either $d$ or $e$, and
$\phi(a) \neq \phi(d), \phi(e)$, so the result follows by the same
argument as in \S \ref{pianotpib}, and the case $n=12$ is done.
}


} 

\bigskip 

{\tiny
Benjamin Howard:
Department of Mathematics,
University of Maryland,
College Park, MD 20742, USA,
bhoward@math.umd.edu

\smallskip

John Millson:
Department of Mathematics,
University of Maryland,
College Park, MD 20742, USA,
jjm@math.umd.edu
 
\smallskip

Andrew Snowden:
Department of Mathematics,
Princeton University,
Princeton, NJ 08544, USA,
asnowden@math.princeton.edu

\smallskip

Ravi Vakil:
Department of Mathematics,
Stanford University, 
Stanford, CA 94305-2125, USA,
vakil@math.stanford.edu
}

\end{document}

%% file: snip.tex
\setlength{\unitlength}{0.00083333in}
\begingroup\makeatletter\ifx\SetFigFont\undefined%
\gdef\SetFigFont#1#2#3#4#5{%
  \reset@font\fontsize{#1}{#2pt}%
  \fontfamily{#3}\fontseries{#4}\fontshape{#5}%
  \selectfont}%
\fi\endgroup%
{\renewcommand{\dashlinestretch}{30}
\begin{picture}(2640,589)(0,-10)
\put(45,529){\blacken\ellipse{74}{74}}
\put(45,529){\ellipse{74}{74}}
\put(45,79){\blacken\ellipse{74}{74}}
\put(45,79){\ellipse{74}{74}}
\put(495,79){\blacken\ellipse{74}{74}}
\put(495,79){\ellipse{74}{74}}
\put(495,529){\blacken\ellipse{74}{74}}
\put(495,529){\ellipse{74}{74}}
\put(1095,529){\blacken\ellipse{74}{74}}
\put(1095,529){\ellipse{74}{74}}
\put(1095,79){\blacken\ellipse{74}{74}}
\put(1095,79){\ellipse{74}{74}}
\put(1545,79){\blacken\ellipse{74}{74}}
\put(1545,79){\ellipse{74}{74}}
\put(1545,529){\blacken\ellipse{74}{74}}
\put(1545,529){\ellipse{74}{74}}
\put(2145,529){\blacken\ellipse{74}{74}}
\put(2145,529){\ellipse{74}{74}}
\put(2145,79){\blacken\ellipse{74}{74}}
\put(2145,79){\ellipse{74}{74}}
\put(2595,79){\blacken\ellipse{74}{74}}
\put(2595,79){\ellipse{74}{74}}
\put(2595,529){\blacken\ellipse{74}{74}}
\put(2595,529){\ellipse{74}{74}}
\path(45,79)(495,79)
\path(45,79)(495,79)
\blacken\path(375.000,49.000)(495.000,79.000)(375.000,109.000)(411.000,79.000)(375.000,49.000)
\path(45,529)(495,529)
\path(45,529)(495,529)
\blacken\path(375.000,499.000)(495.000,529.000)(375.000,559.000)(411.000,529.000)(375.000,499.000)
\path(1095,529)(1095,79)
\blacken\path(1065.000,199.000)(1095.000,79.000)(1125.000,199.000)(1095.000,163.000)(1065.000,199.000)
\path(1545,79)(1095,79)
\blacken\path(1215.000,109.000)(1095.000,79.000)(1215.000,49.000)(1179.000,79.000)(1215.000,109.000)
\path(1095,79)(1545,529)
\blacken\path(1481.360,422.934)(1545.000,529.000)(1438.934,465.360)(1485.603,469.603)(1481.360,422.934)
\path(1545,529)(1545,79)
\blacken\path(1515.000,199.000)(1545.000,79.000)(1575.000,199.000)(1545.000,163.000)(1515.000,199.000)
\path(2145,529)(2595,529)
\blacken\path(2475.000,499.000)(2595.000,529.000)(2475.000,559.000)(2511.000,529.000)(2475.000,499.000)
\path(2145,79)(2595,529)
\blacken\path(2531.360,422.934)(2595.000,529.000)(2488.934,465.360)(2535.603,469.603)(2531.360,422.934)
\path(2145,529)(2145,79)
\blacken\path(2115.000,199.000)(2145.000,79.000)(2175.000,199.000)(2145.000,163.000)(2115.000,199.000)
\path(2595,529)(2595,79)
\blacken\path(2565.000,199.000)(2595.000,79.000)(2625.000,199.000)(2595.000,163.000)(2565.000,199.000)
\path(2595,42)(2145,42)
\blacken\path(2265.000,72.000)(2145.000,42.000)(2265.000,12.000)(2229.000,42.000)(2265.000,72.000)
\path(2145,117)(2595,117)
\blacken\path(2475.000,87.000)(2595.000,117.000)(2475.000,147.000)(2511.000,117.000)(2475.000,87.000)
\put(720,229){\makebox(0,0)[lb]{\smash{{{\SetFigFont{8}{9.6}{\rmdefault}{\mddefault}{\updefault}$\times$}}}}}
\put(1795,229){\makebox(0,0)[lb]{\smash{{{\SetFigFont{8}{9.6}{\rmdefault}{\mddefault}{\updefault}$=$}}}}}
\end{picture}
}

%% file: snoopy.tex
\setlength{\unitlength}{0.00083333in}
\begingroup\makeatletter\ifx\SetFigFont\undefined%
\gdef\SetFigFont#1#2#3#4#5{%
  \reset@font\fontsize{#1}{#2pt}%
  \fontfamily{#3}\fontseries{#4}\fontshape{#5}%
  \selectfont}%
\fi\endgroup%
{\renewcommand{\dashlinestretch}{30}
\begin{picture}(2940,1005)(0,-10)
\put(345,945){\blacken\ellipse{74}{74}}
\put(345,945){\ellipse{74}{74}}
\put(795,945){\blacken\ellipse{74}{74}}
\put(795,945){\ellipse{74}{74}}
\put(1095,495){\blacken\ellipse{74}{74}}
\put(1095,495){\ellipse{74}{74}}
\put(795,45){\blacken\ellipse{74}{74}}
\put(795,45){\ellipse{74}{74}}
\put(345,45){\blacken\ellipse{74}{74}}
\put(345,45){\ellipse{74}{74}}
\put(45,495){\blacken\ellipse{74}{74}}
\put(45,495){\ellipse{74}{74}}
\path(345,945)(795,945)
\path(345,945)(795,945)
\blacken\path(675.000,915.000)(795.000,945.000)(675.000,975.000)(711.000,945.000)(675.000,915.000)
\path(345,45)(45,495)
\path(345,45)(45,495)
\blacken\path(136.526,411.795)(45.000,495.000)(86.603,378.513)(91.595,425.108)(136.526,411.795)
\path(795,45)(1095,495)
\path(795,45)(1095,495)
\blacken\path(1053.397,378.513)(1095.000,495.000)(1003.474,411.795)(1048.405,425.108)(1053.397,378.513)
\put(2145,945){\blacken\ellipse{74}{74}}
\put(2145,945){\ellipse{74}{74}}
\put(2595,945){\blacken\ellipse{74}{74}}
\put(2595,945){\ellipse{74}{74}}
\put(2895,495){\blacken\ellipse{74}{74}}
\put(2895,495){\ellipse{74}{74}}
\put(2595,45){\blacken\ellipse{74}{74}}
\put(2595,45){\ellipse{74}{74}}
\put(2145,45){\blacken\ellipse{74}{74}}
\put(2145,45){\ellipse{74}{74}}
\put(1845,495){\blacken\ellipse{74}{74}}
\put(1845,495){\ellipse{74}{74}}
\path(2145,45)(1845,495)
\path(2145,45)(1845,495)
\blacken\path(1936.526,411.795)(1845.000,495.000)(1886.603,378.513)(1891.595,425.108)(1936.526,411.795)
\path(2595,45)(2895,495)
\path(2595,45)(2895,495)
\blacken\path(2853.397,378.513)(2895.000,495.000)(2803.474,411.795)(2848.405,425.108)(2853.397,378.513)
\blacken\path(2265.000,975.000)(2145.000,945.000)(2265.000,915.000)(2229.000,945.000)(2265.000,975.000)
\path(2145,945)(2595,945)
\path(2145,945)(2595,945)
\put(1245,495){\makebox(0,0)[lb]{\smash{{{\SetFigFont{8}{9.6}{\rmdefault}{\mddefault}{\updefault}$=$}}}}}
\put(1470,495){\makebox(0,0)[lb]{\smash{{{\SetFigFont{8}{9.6}{\rmdefault}{\mddefault}{\updefault}$-$}}}}}
\end{picture}
}

%% file: snide.tex
\setlength{\unitlength}{0.00083333in}
\begingroup\makeatletter\ifx\SetFigFont\undefined%
\gdef\SetFigFont#1#2#3#4#5{%
  \reset@font\fontsize{#1}{#2pt}%
  \fontfamily{#3}\fontseries{#4}\fontshape{#5}%
  \selectfont}%
\fi\endgroup%
{\renewcommand{\dashlinestretch}{30}
\begin{picture}(3000,880)(0,-10)
\put(150,256){\blacken\ellipse{74}{74}}
\put(150,256){\ellipse{74}{74}}
\put(600,706){\blacken\ellipse{74}{74}}
\put(600,706){\ellipse{74}{74}}
\put(1200,256){\blacken\ellipse{74}{74}}
\put(1200,256){\ellipse{74}{74}}
\put(1650,256){\blacken\ellipse{74}{74}}
\put(1650,256){\ellipse{74}{74}}
\put(1650,706){\blacken\ellipse{74}{74}}
\put(1650,706){\ellipse{74}{74}}
\put(2250,256){\blacken\ellipse{74}{74}}
\put(2250,256){\ellipse{74}{74}}
\put(2700,256){\blacken\ellipse{74}{74}}
\put(2700,256){\ellipse{74}{74}}
\put(150,706){\blacken\ellipse{74}{74}}
\put(150,706){\ellipse{74}{74}}
\put(600,256){\blacken\ellipse{74}{74}}
\put(600,256){\ellipse{74}{74}}
\put(1200,706){\blacken\ellipse{74}{74}}
\put(1200,706){\ellipse{74}{74}}
\put(2250,706){\blacken\ellipse{74}{74}}
\put(2250,706){\ellipse{74}{74}}
\put(2700,706){\blacken\ellipse{74}{74}}
\put(2700,706){\ellipse{74}{74}}
\path(150,706)(600,256)
\path(150,706)(600,256)
\blacken\path(493.934,319.640)(600.000,256.000)(536.360,362.066)(540.603,315.397)(493.934,319.640)
\path(1200,706)(1200,256)
\path(1200,706)(1200,256)
\blacken\path(1170.000,376.000)(1200.000,256.000)(1230.000,376.000)(1200.000,340.000)(1170.000,376.000)
\path(2250,706)(2700,706)
\path(2250,706)(2700,706)
\blacken\path(2580.000,676.000)(2700.000,706.000)(2580.000,736.000)(2616.000,706.000)(2580.000,676.000)
\path(2250,256)(2700,256)
\path(2250,256)(2700,256)
\blacken\path(2580.000,226.000)(2700.000,256.000)(2580.000,286.000)(2616.000,256.000)(2580.000,226.000)
\path(600,706)(150,256)
\path(600,706)(150,256)
\blacken\path(213.640,362.066)(150.000,256.000)(256.066,319.640)(209.397,315.397)(213.640,362.066)
\path(1650,256)(1650,706)
\path(1650,256)(1650,706)
\blacken\path(1680.000,586.000)(1650.000,706.000)(1620.000,586.000)(1650.000,622.000)(1680.000,586.000)
\put(675,31){\makebox(0,0)[lb]{\smash{{{\SetFigFont{8}{9.6}{\rmdefault}{\mddefault}{\updefault}$4$}}}}}
\put(0,781){\makebox(0,0)[lb]{\smash{{{\SetFigFont{8}{9.6}{\rmdefault}{\mddefault}{\updefault}$1$}}}}}
\put(2100,781){\makebox(0,0)[lb]{\smash{{{\SetFigFont{8}{9.6}{\rmdefault}{\mddefault}{\updefault}$1$}}}}}
\put(1050,31){\makebox(0,0)[lb]{\smash{{{\SetFigFont{8}{9.6}{\rmdefault}{\mddefault}{\updefault}$3$}}}}}
\put(2100,31){\makebox(0,0)[lb]{\smash{{{\SetFigFont{8}{9.6}{\rmdefault}{\mddefault}{\updefault}$3$}}}}}
\put(1725,31){\makebox(0,0)[lb]{\smash{{{\SetFigFont{8}{9.6}{\rmdefault}{\mddefault}{\updefault}$4$}}}}}
\put(2775,31){\makebox(0,0)[lb]{\smash{{{\SetFigFont{8}{9.6}{\rmdefault}{\mddefault}{\updefault}$4$}}}}}
\put(2775,781){\makebox(0,0)[lb]{\smash{{{\SetFigFont{8}{9.6}{\rmdefault}{\mddefault}{\updefault}$2$}}}}}
\put(675,781){\makebox(0,0)[lb]{\smash{{{\SetFigFont{8}{9.6}{\rmdefault}{\mddefault}{\updefault}$2$}}}}}
\put(0,31){\makebox(0,0)[lb]{\smash{{{\SetFigFont{8}{9.6}{\rmdefault}{\mddefault}{\updefault}$3$}}}}}
\put(1050,781){\makebox(0,0)[lb]{\smash{{{\SetFigFont{8}{9.6}{\rmdefault}{\mddefault}{\updefault}$1$}}}}}
\put(1725,781){\makebox(0,0)[lb]{\smash{{{\SetFigFont{8}{9.6}{\rmdefault}{\mddefault}{\updefault}$2$}}}}}
\put(825,481){\makebox(0,0)[lb]{\smash{{{\SetFigFont{8}{9.6}{\rmdefault}{\mddefault}{\updefault}$+$}}}}}
\put(1875,481){\makebox(0,0)[lb]{\smash{{{\SetFigFont{8}{9.6}{\rmdefault}{\mddefault}{\updefault}$+$}}}}}
\put(3000,481){\makebox(0,0)[lb]{\smash{{{\SetFigFont{8}{9.6}{\rmdefault}{\mddefault}{\updefault}$=0$}}}}}
\end{picture}
}

%% file: snoopdog.tex
\setlength{\unitlength}{0.00083333in}
\begingroup\makeatletter\ifx\SetFigFont\undefined%
\gdef\SetFigFont#1#2#3#4#5{%
  \reset@font\fontsize{#1}{#2pt}%
  \fontfamily{#3}\fontseries{#4}\fontshape{#5}%
  \selectfont}%
\fi\endgroup%
{\renewcommand{\dashlinestretch}{30}
\begin{picture}(4662,999)(0,-10)
\path(3612,942)(4062,942)
\blacken\path(3942.000,912.000)(4062.000,942.000)(3942.000,972.000)(3978.000,942.000)(3942.000,912.000)
\path(12,492)(1062,492)
\blacken\path(942.000,462.000)(1062.000,492.000)(942.000,522.000)(978.000,492.000)(942.000,462.000)
\path(312,942)(762,42)
\blacken\path(681.502,135.915)(762.000,42.000)(735.167,162.748)(724.434,117.132)(681.502,135.915)
\path(762,942)(312,42)
\blacken\path(338.833,162.748)(312.000,42.000)(392.498,135.915)(349.566,117.132)(338.833,162.748)
\path(1662,492)(2712,492)
\blacken\path(2592.000,462.000)(2712.000,492.000)(2592.000,522.000)(2628.000,492.000)(2592.000,462.000)
\path(1962,942)(1962,42)
\blacken\path(1932.000,162.000)(1962.000,42.000)(1992.000,162.000)(1962.000,126.000)(1932.000,162.000)
\path(2412,42)(2412,942)
\blacken\path(2442.000,822.000)(2412.000,942.000)(2382.000,822.000)(2412.000,858.000)(2442.000,822.000)
\path(3612,42)(4062,42)
\blacken\path(3942.000,12.000)(4062.000,42.000)(3942.000,72.000)(3978.000,42.000)(3942.000,12.000)
\path(3312,492)(4362,492)
\blacken\path(4242.000,462.000)(4362.000,492.000)(4242.000,522.000)(4278.000,492.000)(4242.000,462.000)
\put(1287,417){\makebox(0,0)[lb]{\smash{{{\SetFigFont{8}{9.6}{\rmdefault}{\mddefault}{\updefault}$+$}}}}}
\put(2937,417){\makebox(0,0)[lb]{\smash{{{\SetFigFont{8}{9.6}{\rmdefault}{\mddefault}{\updefault}$+$}}}}}
\put(4662,417){\makebox(0,0)[lb]{\smash{{{\SetFigFont{8}{9.6}{\rmdefault}{\mddefault}{\updefault}$=0$}}}}}
\end{picture}
}

%% file: snifter.tex
\setlength{\unitlength}{0.00083333in}
\begingroup\makeatletter\ifx\SetFigFont\undefined%
\gdef\SetFigFont#1#2#3#4#5{%
  \reset@font\fontsize{#1}{#2pt}%
  \fontfamily{#3}\fontseries{#4}\fontshape{#5}%
  \selectfont}%
\fi\endgroup%
{\renewcommand{\dashlinestretch}{30}
\begin{picture}(2400,634)(0,-10)
\put(1725,49){\blacken\ellipse{74}{74}}
\put(1725,49){\ellipse{74}{74}}
\put(2175,49){\blacken\ellipse{74}{74}}
\put(2175,49){\ellipse{74}{74}}
\path(1725,86)(2175,86)
\path(1725,12)(2175,12)
\put(375,574){\blacken\ellipse{74}{74}}
\put(375,574){\ellipse{74}{74}}
\put(375,124){\blacken\ellipse{74}{74}}
\put(375,124){\ellipse{74}{74}}
\put(825,124){\blacken\ellipse{74}{74}}
\put(825,124){\ellipse{74}{74}}
\put(825,574){\blacken\ellipse{74}{74}}
\put(825,574){\ellipse{74}{74}}
\put(1725,574){\blacken\ellipse{74}{74}}
\put(1725,574){\ellipse{74}{74}}
\put(1725,199){\blacken\ellipse{74}{74}}
\put(1725,199){\ellipse{74}{74}}
\put(2175,199){\blacken\ellipse{74}{74}}
\put(2175,199){\ellipse{74}{74}}
\put(2175,574){\blacken\ellipse{74}{74}}
\put(2175,574){\ellipse{74}{74}}
\path(375,574)(825,574)
\path(375,574)(825,574)
\path(375,124)(825,124)
\path(375,124)(825,124)
\path(375,162)(825,162)
\path(375,87)(825,87)
\path(375,574)(375,124)
\path(825,574)(825,124)
\path(1725,199)(1725,574)(2175,574)
	(2175,199)(1725,199)
\path(1575,349)(975,349)
\path(1095.000,379.000)(975.000,349.000)(1095.000,319.000)
\put(0,349){\makebox(0,0)[lb]{\smash{{{\SetFigFont{8}{9.6}{\rmdefault}{\mddefault}{\updefault}$\Gamma$}}}}}
\put(2400,349){\makebox(0,0)[lb]{\smash{{{\SetFigFont{8}{9.6}{\rmdefault}{\mddefault}{\updefault}$\Gamma'$}}}}}
\end{picture}
}

%% file: sniffle.tex
\setlength{\unitlength}{0.00083333in}
\begingroup\makeatletter\ifx\SetFigFont\undefined%
\gdef\SetFigFont#1#2#3#4#5{%
  \reset@font\fontsize{#1}{#2pt}%
  \fontfamily{#3}\fontseries{#4}\fontshape{#5}%
  \selectfont}%
\fi\endgroup%
{\renewcommand{\dashlinestretch}{30}
\begin{picture}(1275,967)(0,-10)
\put(450,781){\blacken\ellipse{74}{74}}
\put(450,781){\ellipse{74}{74}}
\put(900,781){\blacken\ellipse{74}{74}}
\put(900,781){\ellipse{74}{74}}
\put(150,181){\blacken\ellipse{74}{74}}
\put(150,181){\ellipse{74}{74}}
\put(1200,181){\blacken\ellipse{74}{74}}
\put(1200,181){\ellipse{74}{74}}
\path(450,781)(900,781)(150,181)
	(1200,181)(450,781)
\put(852,631){\makebox(0,0)[lb]{\smash{{{\SetFigFont{8}{9.6}{\rmdefault}{\mddefault}{\updefault}$b$}}}}}
\put(975,406){\makebox(0,0)[lb]{\smash{{{\SetFigFont{8}{9.6}{\rmdefault}{\mddefault}{\updefault}$d$}}}}}
\put(300,406){\makebox(0,0)[lb]{\smash{{{\SetFigFont{8}{9.6}{\rmdefault}{\mddefault}{\updefault}$a$}}}}}
\put(637,231){\makebox(0,0)[lb]{\smash{{{\SetFigFont{8}{9.6}{\rmdefault}{\mddefault}{\updefault}$f$}}}}}
\put(637,856){\makebox(0,0)[lb]{\smash{{{\SetFigFont{8}{9.6}{\rmdefault}{\mddefault}{\updefault}$e$}}}}}
\put(400,631){\makebox(0,0)[lb]{\smash{{{\SetFigFont{8}{9.6}{\rmdefault}{\mddefault}{\updefault}$c$}}}}}
\put(975,856){\makebox(0,0)[lb]{\smash{{{\SetFigFont{8}{9.6}{\rmdefault}{\mddefault}{\updefault}$y$}}}}}
\put(1275,31){\makebox(0,0)[lb]{\smash{{{\SetFigFont{8}{9.6}{\rmdefault}{\mddefault}{\updefault}$x$}}}}}
\put(0,31){\makebox(0,0)[lb]{\smash{{{\SetFigFont{8}{9.6}{\rmdefault}{\mddefault}{\updefault}$z$}}}}}
\put(300,856){\makebox(0,0)[lb]{\smash{{{\SetFigFont{8}{9.6}{\rmdefault}{\mddefault}{\updefault}$w$}}}}}
\end{picture}
}

%% file: s5.tex
\setlength{\unitlength}{0.00083333in}
\begingroup\makeatletter\ifx\SetFigFont\undefined%
\gdef\SetFigFont#1#2#3#4#5{%
  \reset@font\fontsize{#1}{#2pt}%
  \fontfamily{#3}\fontseries{#4}\fontshape{#5}%
  \selectfont}%
\fi\endgroup%
{\renewcommand{\dashlinestretch}{30}
\begin{picture}(3750,592)(0,-10)
\path(297,560)(12,353)(121,17)
	(473,17)(582,353)(297,560)
\path(1024,12)(1372,18)(1195,562)(1024,12)
\path(2373,18)(2264,354)(2833,354)
	(2726,21)(2547,565)(2373,18)
\path(3164,349)(3450,559)(3738,349)(3164,349)
\path(3270,15)(3623,15)
\path(3250,77)(3646,77)
\path(1481,317)(909,317)
\path(915,379)(1481,379)
\put(1736,262){\makebox(0,0)[lb]{\smash{{{\SetFigFont{8}{9.6}{\rmdefault}{\mddefault}{\updefault}$=$}}}}}
\end{picture}
}

%% file: s5prep.tex
\setlength{\unitlength}{0.00083333in}
\begingroup\makeatletter\ifx\SetFigFont\undefined%
\gdef\SetFigFont#1#2#3#4#5{%
  \reset@font\fontsize{#1}{#2pt}%
  \fontfamily{#3}\fontseries{#4}\fontshape{#5}%
  \selectfont}%
\fi\endgroup%
{\renewcommand{\dashlinestretch}{30}
\begin{picture}(5844,594)(0,-10)
\path(467,351)(751,567)(1034,357)
\path(574,23)(927,23)
\path(5268,346)(5832,349)
\path(5377,12)(5719,15)
\path(3765,349)(3875,15)
\path(4225,12)(4332,352)
\path(1965,346)(2537,346)
\path(2077,12)(2251,559)(2420,15)
\put(0,259){\makebox(0,0)[lb]{\smash{{{\SetFigFont{8}{9.6}{\rmdefault}{\mddefault}{\updefault}$\Gamma_1=$}}}}}
\put(1350,259){\makebox(0,0)[lb]{\smash{{{\SetFigFont{8}{9.6}{\rmdefault}{\mddefault}{\updefault}$\Gamma_2=$}}}}}
\put(3150,259){\makebox(0,0)[lb]{\smash{{{\SetFigFont{8}{9.6}{\rmdefault}{\mddefault}{\updefault}$\Delta_1=$}}}}}
\put(4650,259){\makebox(0,0)[lb]{\smash{{{\SetFigFont{8}{9.6}{\rmdefault}{\mddefault}{\updefault}$\Delta_2=$}}}}}
\end{picture}
}

%% file: s835.tex
\setlength{\unitlength}{0.00083333in}
\begingroup\makeatletter\ifx\SetFigFont\undefined%
\gdef\SetFigFont#1#2#3#4#5{%
  \reset@font\fontsize{#1}{#2pt}%
  \fontfamily{#3}\fontseries{#4}\fontshape{#5}%
  \selectfont}%
\fi\endgroup%
{\renewcommand{\dashlinestretch}{30}
\begin{picture}(3636,649)(0,-10)
\path(189,618)(12,443)
\path(441,619)(615,442)
\path(619,193)(440,17)
\path(191,16)(12,198)
\path(1093,622)(1346,622)
\path(915,447)(1517,448)
\path(914,196)(1520,196)
\path(1093,20)(1346,20)
\path(2300,618)(2123,442)
\path(2550,621)(2731,439)
\path(2124,196)(2726,196)
\path(2300,17)(2554,17)
\path(3194,618)(3446,618)
\path(3022,442)(3622,442)
\path(3019,195)(3197,16)
\path(3443,12)(3624,195)
\put(1714,301){\makebox(0,0)[lb]{\smash{{{\SetFigFont{8}{9.6}{\rmdefault}{\mddefault}{\updefault}$=$}}}}}
\end{picture}
}

%% file: snegre.tex
\setlength{\unitlength}{0.00083333in}
\begingroup\makeatletter\ifx\SetFigFont\undefined%
\gdef\SetFigFont#1#2#3#4#5{%
  \reset@font\fontsize{#1}{#2pt}%
  \fontfamily{#3}\fontseries{#4}\fontshape{#5}%
  \selectfont}%
\fi\endgroup%
{\renewcommand{\dashlinestretch}{30}
\begin{picture}(5419,566)(0,-10)
\path(161,534)(457,534)
\path(12,271)(607,274)
\path(161,12)(457,18)
\path(911,274)(1058,18)
\path(1058,534)(1357,18)
\path(1360,534)(1506,280)
\path(1810,274)(1959,539)
\path(1959,15)(2258,539)
\path(2258,15)(2407,277)
\path(3158,534)(3457,18)
\path(3158,15)(3457,539)
\path(3009,280)(3609,280)
\path(4061,534)(4359,536)
\path(3908,277)(4058,18)
\path(4359,18)(4509,277)
\path(4959,531)(4810,280)
\path(4959,15)(5258,15)
\path(5261,534)(5407,277)
\put(2631,243){\makebox(0,0)[lb]{\smash{{{\SetFigFont{8}{9.6}{\rmdefault}{\mddefault}{\updefault}$=$}}}}}
\end{picture}
}

%% file: snegreplus.tex
\setlength{\unitlength}{0.00083333in}
\begingroup\makeatletter\ifx\SetFigFont\undefined%
\gdef\SetFigFont#1#2#3#4#5{%
  \reset@font\fontsize{#1}{#2pt}%
  \fontfamily{#3}\fontseries{#4}\fontshape{#5}%
  \selectfont}%
\fi\endgroup%
{\renewcommand{\dashlinestretch}{30}
\begin{picture}(5419,751)(0,-10)
\path(161,719)(457,719)
\path(12,456)(607,459)
\path(161,197)(457,203)
\path(911,459)(1058,203)
\path(1058,719)(1357,203)
\path(1360,719)(1506,465)
\path(1810,459)(1959,724)
\path(1959,200)(2258,724)
\path(2258,200)(2407,462)
\path(3158,719)(3457,203)
\path(3158,200)(3457,724)
\path(3009,465)(3609,465)
\path(4061,719)(4359,721)
\path(3908,462)(4058,203)
\path(4359,203)(4509,462)
\path(4959,716)(4810,465)
\path(4959,200)(5258,200)
\path(5261,719)(5407,462)
\path(158,12)(458,12)
\path(1058,12)(1358,12)
\path(1958,12)(2258,12)
\path(3158,12)(3458,12)
\path(4058,12)(4358,12)
\path(4958,12)(5258,12)
\put(2618,338){\makebox(0,0)[lb]{\smash{{{\SetFigFont{8}{9.6}{\rmdefault}{\mddefault}{\updefault}$=$}}}}}
\end{picture}
}

%% file: m6.tex
\setlength{\unitlength}{0.00083333in}
\begingroup\makeatletter\ifx\SetFigFont\undefined%
\gdef\SetFigFont#1#2#3#4#5{%
  \reset@font\fontsize{#1}{#2pt}%
  \fontfamily{#3}\fontseries{#4}\fontshape{#5}%
  \selectfont}%
\fi\endgroup%
{\renewcommand{\dashlinestretch}{30}
\begin{picture}(2850,1675)(0,-10)
\path(1125,810)(1275,919)(1181,636)
\path(1367,635)(1423,812)
\put(1236,964){\makebox(0,0)[lb]{\smash{{{\SetFigFont{8}{9.6}{\rmdefault}{\mddefault}{\updefault}$2$}}}}}
\path(1936,1648)(2089,1648)
\path(1862,1519)(1936,1391)
\path(2085,1390)(2160,1519)
\put(1492,1538){\makebox(0,0)[lb]{\smash{{{\SetFigFont{8}{9.6}{\rmdefault}{\mddefault}{\updefault}$\Gamma =$}}}}}
\put(2542,1538){\makebox(0,0)[lb]{\smash{{{\SetFigFont{8}{9.6}{\rmdefault}{\mddefault}{\updefault}$\operatorname{deg}=3$}}}}}
\path(225,1461)(225,936)
\blacken\path(195.000,1056.000)(225.000,936.000)(255.000,1056.000)(225.000,1020.000)(195.000,1056.000)
\path(225,711)(225,186)
\blacken\path(195.000,306.000)(225.000,186.000)(255.000,306.000)(225.000,270.000)(195.000,306.000)
\path(686,76)(1725,76)
\blacken\path(1605.000,46.000)(1725.000,76.000)(1605.000,106.000)(1641.000,76.000)(1605.000,46.000)
\put(75,1536){\makebox(0,0)[lb]{\smash{{{\SetFigFont{8}{9.6}{\rmdefault}{\mddefault}{\updefault}${\mathbf{w}}=(1,1,1,1,1,1)$}}}}}
\put(0,1161){\makebox(0,0)[lb]{\smash{{{\SetFigFont{8}{9.6}{\rmdefault}{\mddefault}{\updefault}(c)}}}}}
\put(375,1161){\makebox(0,0)[lb]{\smash{{{\SetFigFont{8}{9.6}{\rmdefault}{\mddefault}{\updefault}$\times 3$}}}}}
\put(75,786){\makebox(0,0)[lb]{\smash{{{\SetFigFont{8}{9.6}{\rmdefault}{\mddefault}{\updefault}$(2,1,1,1,1)$}}}}}
\put(75,36){\makebox(0,0)[lb]{\smash{{{\SetFigFont{8}{9.6}{\rmdefault}{\mddefault}{\updefault}$(2,2,1,1)$}}}}}
\put(0,411){\makebox(0,0)[lb]{\smash{{{\SetFigFont{8}{9.6}{\rmdefault}{\mddefault}{\updefault}(c)}}}}}
\put(1800,36){\makebox(0,0)[lb]{\smash{{{\SetFigFont{8}{9.6}{\rmdefault}{\mddefault}{\updefault}$(1,1,1,1)$}}}}}
\put(1725,786){\makebox(0,0)[lb]{\smash{{{\SetFigFont{8}{9.6}{\rmdefault}{\mddefault}{\updefault}$\operatorname{deg}=1$}}}}}
\put(2700,36){\makebox(0,0)[lb]{\smash{{{\SetFigFont{8}{9.6}{\rmdefault}{\mddefault}{\updefault}$\operatorname{deg}=1$}}}}}
\put(1050,186){\makebox(0,0)[lb]{\smash{{{\SetFigFont{8}{9.6}{\rmdefault}{\mddefault}{\updefault}(d)}}}}}
\put(2850,186){\makebox(0,0)[lb]{\smash{{{\SetFigFont{8}{9.6}{\rmdefault}{\mddefault}{\updefault}(b)}}}}}
\end{picture}
}

%% file: m5.tex
\setlength{\unitlength}{0.00083333in}
\begingroup\makeatletter\ifx\SetFigFont\undefined%
\gdef\SetFigFont#1#2#3#4#5{%
  \reset@font\fontsize{#1}{#2pt}%
  \fontfamily{#3}\fontseries{#4}\fontshape{#5}%
  \selectfont}%
\fi\endgroup%
{\renewcommand{\dashlinestretch}{30}
\begin{picture}(4575,973)(0,-10)
\path(1730,819)(2079,817)(1904,946)(1730,819)
\path(1792,614)(2013,614)
\path(1785,642)(2022,642)
\put(1365,767){\makebox(0,0)[lb]{\smash{{{\SetFigFont{8}{9.6}{\rmdefault}{\mddefault}{\updefault}$\Gamma =$}}}}}
\put(2415,767){\makebox(0,0)[lb]{\smash{{{\SetFigFont{8}{9.6}{\rmdefault}{\mddefault}{\updefault}$\operatorname{deg}=5$}}}}}
\path(225,711)(225,186)
\blacken\path(195.000,306.000)(225.000,186.000)(255.000,306.000)(225.000,270.000)(195.000,306.000)
\path(675,111)(1200,111)
\blacken\path(1080.000,81.000)(1200.000,111.000)(1080.000,141.000)(1116.000,111.000)(1080.000,81.000)
\path(1875,111)(2475,111)
\blacken\path(2355.000,81.000)(2475.000,111.000)(2355.000,141.000)(2391.000,111.000)(2355.000,81.000)
\path(3150,111)(3600,111)
\blacken\path(3480.000,81.000)(3600.000,111.000)(3480.000,141.000)(3516.000,111.000)(3480.000,81.000)
\put(0,411){\makebox(0,0)[lb]{\smash{{{\SetFigFont{8}{9.6}{\rmdefault}{\mddefault}{\updefault}(c)}}}}}
\put(75,36){\makebox(0,0)[lb]{\smash{{{\SetFigFont{8}{9.6}{\rmdefault}{\mddefault}{\updefault}$(4,2,2,2)$}}}}}
\put(75,786){\makebox(0,0)[lb]{\smash{{{\SetFigFont{8}{9.6}{\rmdefault}{\mddefault}{\updefault}${\mathbf{w}}=(2,2,2,2,2)$}}}}}
\put(375,411){\makebox(0,0)[lb]{\smash{{{\SetFigFont{8}{9.6}{\rmdefault}{\mddefault}{\updefault}$\times 5$}}}}}
\put(1275,36){\makebox(0,0)[lb]{\smash{{{\SetFigFont{8}{9.6}{\rmdefault}{\mddefault}{\updefault}$(3,2,2,1)$}}}}}
\put(2550,36){\makebox(0,0)[lb]{\smash{{{\SetFigFont{8}{9.6}{\rmdefault}{\mddefault}{\updefault}$(2,2,1,1)$}}}}}
\put(3675,36){\makebox(0,0)[lb]{\smash{{{\SetFigFont{8}{9.6}{\rmdefault}{\mddefault}{\updefault}$(1,1,1,1)$}}}}}
\put(4500,36){\makebox(0,0)[lb]{\smash{{{\SetFigFont{8}{9.6}{\rmdefault}{\mddefault}{\updefault}$\operatorname{deg}=1$}}}}}
\put(825,186){\makebox(0,0)[lb]{\smash{{{\SetFigFont{8}{9.6}{\rmdefault}{\mddefault}{\updefault}(d)}}}}}
\put(2100,186){\makebox(0,0)[lb]{\smash{{{\SetFigFont{8}{9.6}{\rmdefault}{\mddefault}{\updefault}(d)}}}}}
\put(3225,186){\makebox(0,0)[lb]{\smash{{{\SetFigFont{8}{9.6}{\rmdefault}{\mddefault}{\updefault}(d)}}}}}
\put(4575,186){\makebox(0,0)[lb]{\smash{{{\SetFigFont{8}{9.6}{\rmdefault}{\mddefault}{\updefault}(b)}}}}}
\end{picture}
}

%% file: snape.tex
\setlength{\unitlength}{0.00083333in}
\begingroup\makeatletter\ifx\SetFigFont\undefined%
\gdef\SetFigFont#1#2#3#4#5{%
  \reset@font\fontsize{#1}{#2pt}%
  \fontfamily{#3}\fontseries{#4}\fontshape{#5}%
  \selectfont}%
\fi\endgroup%
{\renewcommand{\dashlinestretch}{30}
\begin{picture}(1902,1044)(0,-10)
\put(975,795){\blacken\ellipse{74}{74}}
\put(975,795){\ellipse{74}{74}}
\put(375,645){\blacken\ellipse{74}{74}}
\put(375,645){\ellipse{74}{74}}
\put(1575,645){\blacken\ellipse{74}{74}}
\put(1575,645){\ellipse{74}{74}}
\put(225,45){\blacken\ellipse{74}{74}}
\put(225,45){\ellipse{74}{74}}
\put(1725,45){\blacken\ellipse{74}{74}}
\put(1725,45){\ellipse{74}{74}}
\path(225,45)(375,645)(975,795)
	(1575,645)(1725,45)
\put(0,45){\makebox(0,0)[lb]{\smash{{{\SetFigFont{8}{9.6}{\rmdefault}{\mddefault}{\updefault}$p$}}}}}
\put(150,720){\makebox(0,0)[lb]{\smash{{{\SetFigFont{8}{9.6}{\rmdefault}{\mddefault}{\updefault}$q$}}}}}
\put(1725,720){\makebox(0,0)[lb]{\smash{{{\SetFigFont{8}{9.6}{\rmdefault}{\mddefault}{\updefault}$s$}}}}}
\put(125,346){\makebox(0,0)[lb]{\smash{{{\SetFigFont{8}{9.6}{\rmdefault}{\mddefault}{\updefault}$\Delta$}}}}}
\put(1688,346){\makebox(0,0)[lb]{\smash{{{\SetFigFont{8}{9.6}{\rmdefault}{\mddefault}{\updefault}$\Gamma$}}}}}
\put(1902,47){\makebox(0,0)[lb]{\smash{{{\SetFigFont{8}{9.6}{\rmdefault}{\mddefault}{\updefault}$t$}}}}}
\put(927,945){\makebox(0,0)[lb]{\smash{{{\SetFigFont{8}{9.6}{\rmdefault}{\mddefault}{\updefault}$r$}}}}}
\put(1246,772){\makebox(0,0)[lb]{\smash{{{\SetFigFont{8}{9.6}{\rmdefault}{\mddefault}{\updefault}$\Delta$}}}}}
\put(583,772){\makebox(0,0)[lb]{\smash{{{\SetFigFont{8}{9.6}{\rmdefault}{\mddefault}{\updefault}$\Gamma$}}}}}
\end{picture}
}

%% file: 10case.tex
\setlength{\unitlength}{0.00083333in}
\begingroup\makeatletter\ifx\SetFigFont\undefined%
\gdef\SetFigFont#1#2#3#4#5{%
  \reset@font\fontsize{#1}{#2pt}%
  \fontfamily{#3}\fontseries{#4}\fontshape{#5}%
  \selectfont}%
\fi\endgroup%
{\renewcommand{\dashlinestretch}{30}
\begin{picture}(3128,967)(0,-10)
\put(83,181){\blacken\ellipse{74}{74}}
\put(83,181){\ellipse{74}{74}}
\put(683,781){\blacken\ellipse{74}{74}}
\put(683,781){\ellipse{74}{74}}
\put(683,181){\blacken\ellipse{74}{74}}
\put(683,181){\ellipse{74}{74}}
\put(1283,181){\blacken\ellipse{74}{74}}
\put(1283,181){\ellipse{74}{74}}
\put(1883,781){\blacken\ellipse{74}{74}}
\put(1883,781){\ellipse{74}{74}}
\put(1883,181){\blacken\ellipse{74}{74}}
\put(1883,181){\ellipse{74}{74}}
\put(2483,181){\blacken\ellipse{74}{74}}
\put(2483,181){\ellipse{74}{74}}
\put(3083,181){\blacken\ellipse{74}{74}}
\put(3083,181){\ellipse{74}{74}}
\put(83,781){\blacken\ellipse{74}{74}}
\put(83,781){\ellipse{74}{74}}
\put(1283,781){\blacken\ellipse{74}{74}}
\put(1283,781){\ellipse{74}{74}}
\path(83,181)(83,781)(683,781)(683,181)
\path(1283,181)(1283,781)(1883,781)(1883,181)
\path(2483,181)(3083,181)
\put(683,31){\makebox(0,0)[lb]{\smash{{{\SetFigFont{8}{9.6}{\rmdefault}{\mddefault}{\updefault}$b$}}}}}
\put(1883,31){\makebox(0,0)[lb]{\smash{{{\SetFigFont{8}{9.6}{\rmdefault}{\mddefault}{\updefault}$d$}}}}}
\put(3083,31){\makebox(0,0)[lb]{\smash{{{\SetFigFont{8}{9.6}{\rmdefault}{\mddefault}{\updefault}$f$}}}}}
\put(1883,856){\makebox(0,0)[lb]{\smash{{{\SetFigFont{8}{9.6}{\rmdefault}{\mddefault}{\updefault}$j$}}}}}
\put(2400,31){\makebox(0,0)[lb]{\smash{{{\SetFigFont{8}{9.6}{\rmdefault}{\mddefault}{\updefault}$e$}}}}}
\put(1200,31){\makebox(0,0)[lb]{\smash{{{\SetFigFont{8}{9.6}{\rmdefault}{\mddefault}{\updefault}$c$}}}}}
\put(0,31){\makebox(0,0)[lb]{\smash{{{\SetFigFont{8}{9.6}{\rmdefault}{\mddefault}{\updefault}$a$}}}}}
\put(2725,213){\makebox(0,0)[lb]{\smash{{{\SetFigFont{8}{9.6}{\rmdefault}{\mddefault}{\updefault}$\Gamma$}}}}}
\put(1737,444){\makebox(0,0)[lb]{\smash{{{\SetFigFont{8}{9.6}{\rmdefault}{\mddefault}{\updefault}$\Gamma$}}}}}
\put(561,444){\makebox(0,0)[lb]{\smash{{{\SetFigFont{8}{9.6}{\rmdefault}{\mddefault}{\updefault}$\Gamma$}}}}}
\put(135,444){\makebox(0,0)[lb]{\smash{{{\SetFigFont{8}{9.6}{\rmdefault}{\mddefault}{\updefault}$\Gamma$}}}}}
\put(320,839){\makebox(0,0)[lb]{\smash{{{\SetFigFont{8}{9.6}{\rmdefault}{\mddefault}{\updefault}$\Delta$}}}}}
\put(1522,841){\makebox(0,0)[lb]{\smash{{{\SetFigFont{8}{9.6}{\rmdefault}{\mddefault}{\updefault}$\Delta$}}}}}
\put(1337,444){\makebox(0,0)[lb]{\smash{{{\SetFigFont{8}{9.6}{\rmdefault}{\mddefault}{\updefault}$\Gamma$}}}}}
\put(0,856){\makebox(0,0)[lb]{\smash{{{\SetFigFont{8}{9.6}{\rmdefault}{\mddefault}{\updefault}$g$}}}}}
\put(683,856){\makebox(0,0)[lb]{\smash{{{\SetFigFont{8}{9.6}{\rmdefault}{\mddefault}{\updefault}$h$}}}}}
\put(1200,856){\makebox(0,0)[lb]{\smash{{{\SetFigFont{8}{9.6}{\rmdefault}{\mddefault}{\updefault}$i$}}}}}
\end{picture}
}

%% file: 12case.tex
\setlength{\unitlength}{0.00083333in}
\begingroup\makeatletter\ifx\SetFigFont\undefined%
\gdef\SetFigFont#1#2#3#4#5{%
  \reset@font\fontsize{#1}{#2pt}%
  \fontfamily{#3}\fontseries{#4}\fontshape{#5}%
  \selectfont}%
\fi\endgroup%
{\renewcommand{\dashlinestretch}{30}
\begin{picture}(3128,967)(0,-10)
\put(83,181){\blacken\ellipse{74}{74}}
\put(83,181){\ellipse{74}{74}}
\put(83,781){\blacken\ellipse{74}{74}}
\put(83,781){\ellipse{74}{74}}
\put(683,781){\blacken\ellipse{74}{74}}
\put(683,781){\ellipse{74}{74}}
\put(683,181){\blacken\ellipse{74}{74}}
\put(683,181){\ellipse{74}{74}}
\put(1283,181){\blacken\ellipse{74}{74}}
\put(1283,181){\ellipse{74}{74}}
\put(1283,781){\blacken\ellipse{74}{74}}
\put(1283,781){\ellipse{74}{74}}
\put(1883,781){\blacken\ellipse{74}{74}}
\put(1883,781){\ellipse{74}{74}}
\put(1883,181){\blacken\ellipse{74}{74}}
\put(1883,181){\ellipse{74}{74}}
\put(2483,181){\blacken\ellipse{74}{74}}
\put(2483,181){\ellipse{74}{74}}
\put(3083,181){\blacken\ellipse{74}{74}}
\put(3083,181){\ellipse{74}{74}}
\put(2483,781){\blacken\ellipse{74}{74}}
\put(2483,781){\ellipse{74}{74}}
\put(3083,781){\blacken\ellipse{74}{74}}
\put(3083,781){\ellipse{74}{74}}
\path(83,181)(83,781)(683,781)(683,181)
\path(1283,181)(1283,781)(1883,781)(1883,181)
\path(2483,181)(2483,781)(3083,781)(3083,181)
\put(0,856){\makebox(0,0)[lb]{\smash{{{\SetFigFont{8}{9.6}{\rmdefault}{\mddefault}{\updefault}$g$}}}}}
\put(683,31){\makebox(0,0)[lb]{\smash{{{\SetFigFont{8}{9.6}{\rmdefault}{\mddefault}{\updefault}$b$}}}}}
\put(1883,31){\makebox(0,0)[lb]{\smash{{{\SetFigFont{8}{9.6}{\rmdefault}{\mddefault}{\updefault}$d$}}}}}
\put(3083,31){\makebox(0,0)[lb]{\smash{{{\SetFigFont{8}{9.6}{\rmdefault}{\mddefault}{\updefault}$f$}}}}}
\put(1200,856){\makebox(0,0)[lb]{\smash{{{\SetFigFont{8}{9.6}{\rmdefault}{\mddefault}{\updefault}$i$}}}}}
\put(1883,856){\makebox(0,0)[lb]{\smash{{{\SetFigFont{8}{9.6}{\rmdefault}{\mddefault}{\updefault}$j$}}}}}
\put(2400,31){\makebox(0,0)[lb]{\smash{{{\SetFigFont{8}{9.6}{\rmdefault}{\mddefault}{\updefault}$e$}}}}}
\put(1200,31){\makebox(0,0)[lb]{\smash{{{\SetFigFont{8}{9.6}{\rmdefault}{\mddefault}{\updefault}$c$}}}}}
\put(0,31){\makebox(0,0)[lb]{\smash{{{\SetFigFont{8}{9.6}{\rmdefault}{\mddefault}{\updefault}$a$}}}}}
\put(1337,444){\makebox(0,0)[lb]{\smash{{{\SetFigFont{8}{9.6}{\rmdefault}{\mddefault}{\updefault}$\Gamma$}}}}}
\put(1737,444){\makebox(0,0)[lb]{\smash{{{\SetFigFont{8}{9.6}{\rmdefault}{\mddefault}{\updefault}$\Gamma$}}}}}
\put(561,444){\makebox(0,0)[lb]{\smash{{{\SetFigFont{8}{9.6}{\rmdefault}{\mddefault}{\updefault}$\Gamma$}}}}}
\put(135,444){\makebox(0,0)[lb]{\smash{{{\SetFigFont{8}{9.6}{\rmdefault}{\mddefault}{\updefault}$\Gamma$}}}}}
\put(320,839){\makebox(0,0)[lb]{\smash{{{\SetFigFont{8}{9.6}{\rmdefault}{\mddefault}{\updefault}$\Delta$}}}}}
\put(2531,444){\makebox(0,0)[lb]{\smash{{{\SetFigFont{8}{9.6}{\rmdefault}{\mddefault}{\updefault}$\Gamma$}}}}}
\put(2945,444){\makebox(0,0)[lb]{\smash{{{\SetFigFont{8}{9.6}{\rmdefault}{\mddefault}{\updefault}$\Gamma$}}}}}
\put(2719,837){\makebox(0,0)[lb]{\smash{{{\SetFigFont{8}{9.6}{\rmdefault}{\mddefault}{\updefault}$\Delta$}}}}}
\put(2400,856){\makebox(0,0)[lb]{\smash{{{\SetFigFont{8}{9.6}{\rmdefault}{\mddefault}{\updefault}$k$}}}}}
\put(683,856){\makebox(0,0)[lb]{\smash{{{\SetFigFont{8}{9.6}{\rmdefault}{\mddefault}{\updefault}$h$}}}}}
\put(1522,841){\makebox(0,0)[lb]{\smash{{{\SetFigFont{8}{9.6}{\rmdefault}{\mddefault}{\updefault}$\Delta$}}}}}
\put(3083,856){\makebox(0,0)[lb]{\smash{{{\SetFigFont{8}{9.6}{\rmdefault}{\mddefault}{\updefault}$l$}}}}}
\end{picture}
}